# ON ASYMPTOTICALLY OPTIMAL TESTS UNDER LOSS OF IDENTIFIABILITY IN SEMIPARAMETRIC MODELS

By Rui Song[1], Michael R. Kosorok[1] and Jason P. Fine[2]

*University of North Carolina*

We consider tests of hypotheses when the parameters are not identifiable under the null in semiparametric models, where regularity conditions for profile likelihood theory fail. Exponential average tests based on integrated profile likelihood are constructed and shown to be asymptotically optimal under a weighted average power criterion with respect to a prior on the nonidentifiable aspect of the model. These results extend existing results for parametric models, which involve more restrictive assumptions on the form of the alternative than do our results. Moreover, the proposed tests accommodate models with infinite dimensional nuisance parameters which either may not be identifiable or may not be estimable at the usual parametric rate. Examples include tests of the presence of a change-point in the Cox model with current status data and tests of regression parameters in odds-rate models with right censored data. Optimal tests have not previously been studied for these scenarios. We study the asymptotic distribution of the proposed tests under the null, fixed contiguous alternatives and random contiguous alternatives. We also propose a weighted bootstrap procedure for computing the critical values of the test statistics. The optimal tests perform well in simulation studies, where they may exhibit improved power over alternative tests.

**1. Introduction.** In this paper we investigate nonstandard testing problems involving a family of probability distributions $\{P_\theta, \theta \in \Theta\}$, known up to a parameter $\theta$, in a parameter space $\Theta$. The parameter space $\Theta$ is assumed to be a subset of an infinite-dimensional metric space. The null and alternative hypotheses are:

$$H_0 : \theta \in \Theta_0 \quad \text{vs.} \quad H_1 : \theta \in \Theta \backslash \Theta_0,$$

Received October 2007; revised June 2008.

[1]Supported in part by Grant CA075142 from the National Cancer Institute.

[2]Supported in part by Grant CA094893 from the National Cancer Institute.

*AMS 2000 subject classifications.* Primary 62A01, 62G10; secondary 62G20, 62C99.

*Key words and phrases.* Change-point models, contiguous alternative, empirical processes, exponential average test, nonstandard testing problem, odds-rate models, optimal test, power, profile likelihood.

where $\Theta_0$ is a subset of $\Theta$ and contains at least two elements. In the usual testing framework, the parameters are unique under the null so that identifiability is not an issue. While we allow multiple values of $\theta$ satisfying the null, we assume that the null distribution, denoted by $P_0$, is unique, where $\Theta_0 = \{\theta \in \Theta : P_\theta = P_0\}$. Under this setup, the true value of $\theta$ is not identifiable under the null, since for any $\theta \neq \theta'$ in $\Theta_0$, $P_\theta = P_{\theta'} = P_0$. Such loss of identifiability occurs in diverse applications in the social, biological, physical and medical sciences. We next present two such examples followed by a description of the main contributions of this paper. The Introduction concludes with a brief outline of the remainder of the paper.

1.1. *Example 1: Univariate frailty regression under right censoring.* Let $T$ be a nonnegative random variable representing the failure time, $C$ be the independent censoring time, $V \equiv \min(T, C)$ and $Z \equiv Z(\cdot)$ be a corresponding $p$-dimensional covariate process. The observed data $\{X_i = (V_i, \Delta_i, Z_i), i = 1, \ldots, n\}$ consists of $n$ i.i.d. realizations of $X = (V, \Delta, Z)$, where $\Delta \equiv 1\{T \leq V\}$, $1\{\cdot\}$ is the indicator function. In this model, the hazard function of the survival time $T$ given covariates $Z$ is

$$(1) \qquad \lambda\{t; Z(t), W\} = \eta(t) W \exp\{\beta^T Z(t)\},$$

where $t$ is the time index, $W$ is an unobserved gamma frailty with mean 1 and variance $\zeta$, $\beta$ is a $p$-dimensional regression parameter and $\eta(\cdot)$ is a completely unspecified baseline hazard function.

When $\beta$ is not zero, the odds-rate model has been treated extensively; see Kosorok, Lee and Fine (2004), Murphy, Rossini and van der Vaart (1997); Murphy and van der Vaart (1997, 2000); Parner (1998); Slud and Vonta (2004), among others. Scharfstein, Tsiatis and Gilbert (1998) considered semi-parametric efficient estimation in the setting, where the covariates are time independent, $\zeta$ is assumed known and $\eta(\cdot)$ is assumed to be absolutely continuous. Bagdonavičius and Nikulin (1999) considered estimation for a class of proportional hazards model, which includes the odds-rate model with $\zeta$ unspecified, based on a modified partial likelihood. Kosorok, Lee and Fine (2004) considered robust inference for odds-rate models when the frailty distribution and regression covariates may be misspecified. To our knowledge, problems associated with testing the null $\beta = 0$ when the frailty parameter is unknown have not been previously considered in the statistical literature.

It has been shown that $\zeta$ and $\eta(\cdot)$ are not identifiable under the null [Kosorok, Lee and Fine (2004)]. Intuitively, when $\beta = 0$, the covariate process $Z$ provides no information for the failure time process. The frailty $W$ and the baseline hazard $\eta(\cdot)$ are not distinguishable from each other, hence $\zeta$ and $\eta(\cdot)$ are not identifiable. Thus, the testing problem described above is nonregular and standard asymptotic results are not applicable.

1.2. *Example 2: Change-point regression for current status data.* Change-point models have been studied extensively and have proven to be popular in clinical research. In many settings, a change-point effect is realistic and can be much easier to interpret than a quadratic or more complex nonlinear effect [Chappell (1989)]. Change-point Cox models have been widely used in survival applications, as in Kosorok and Song (2007); Luo, Turnbull and Clark (1997); Pons (2003), where likelihood ratio tests were investigated. However, to our knowledge, optimal testing has not been explored for such models.

Under current status censoring, a subject is examined once at a random observation time $V$ and at that time it is observed whether the event time $T \leq V$ or not. The observed data $\{X_i = (V_i, \Delta_i, Z_i), i = 1, \ldots, n\}$ consists of $n$ i.i.d. realizations of $X = (V, \Delta, Z)$, where $\Delta \equiv 1\{T \leq V\}$ and $Z$ is a $d$-dimensional covariate. Here we let $d = 1$ for simplicity. In this example, we assume that the time to event $T$ satisfies a change-point Cox model conditionally on the covariate $Z$. That is, the density of $X$ is given by:

$$(2) \qquad p_\theta(x) = (1 - e^{-e^{r_\gamma(z)}\Lambda(v)})^\Delta (e^{-e^{r_\gamma(z)}\Lambda(v)})^{1-\Delta} f_{V,Z}(v, z),$$

with $r_\gamma(z) = \alpha z + (\beta_1 + \beta_2 z)1\{z > \zeta\}$, where $\alpha$, $\beta_1$ and $\beta_2$ are scalar regression parameters, $\zeta$ is the change-point parameter and $\Lambda(\cdot)$ is the cumulative baseline hazard function. We also define the collected parameters $\beta \equiv (\beta_1, \beta_2)$, $\xi \equiv (\beta, \alpha)$, $\gamma \equiv (\xi, \zeta)$ and $\eta \equiv (\alpha, \Lambda)$. We are particularly interested in the hypothesis test of the existence of a change-point for regression parameters in the above model, that is, $H_0 : \beta = 0$.

Although Cox regression with current status data was discussed by Huang (1996) and others, change-point Cox regression has not been studied with current status data. The development of optimal tests in the current status setting is further complicated by the fact that the nuisance parameter $\Lambda$ cannot be estimated at the parametric rate, unlike with right censored data.

In model (2), the change-point parameter is present only under the alternative. This is different from Example 1, where the odds rate parameter $\zeta$ and the baseline hazard function $\eta(\cdot)$ are both present, but indistinguishable, under the null.

1.3. *Description of main contribution.* The statistical literature contains numerous precedents on the nonidentifiability problem in parametric models, see Chernoff (1954), Chernoff and Lander (1995), Dacunha-Castelle and Gassiat (1999) and Liu and Shao (2003). Among others, Dacunha-Castelle and Gassiat (1999) proposed a locally conic parametrization approach to enable asymptotic expansions of the likelihood ratio test under loss of identifiability under the null. Liu and Shao (2003) derived a quadratic approximation of the log-likelihood ratio function by using Hellinger distance. Most authors directly study the approximation of the log-likelihood ratio function in some neighborhood and obtain its asymptotic null distribution. However, the asymptotic optimality properties of the classical likelihood ratio tests (LRT) do not hold anymore [Lindsay (1995)] and Wald and score tests are not even well defined in these nonstandard problems. To our knowledge, all results for testing nonidentifiable $P_0$ using likelihood based tests are for parametric models. The main aim of this paper is to investigate the construction of optimal likelihood based tests for semiparametric models.

A key question which arises, as noted by Dacunha-Castelle and Gassiat (1999), is: since the parameter is not identifiable, around which point can an expansion be made? To address this question, we assume the existence of a "full rank" reparameterization which contains all the information of the null model and in which all parameters are identifiable. To be specific, we partition $\theta \equiv (\psi, \zeta)$ and $\psi \equiv (\beta, \eta)$, where $\beta \in \mathbb{R}^p$ is a parameter of interest, $\zeta \in \mathbb{R}^q$ and $\eta$ is a parameter defined on an arbitrary parametric space, $\mathcal{H}_\eta$. We assume that the information in the null model can be absorbed into the parameter space of $\eta$, through this full rank reparameterization. This is made precise in Section 2. Note that Example 1 requires such a reparameterization since both $\zeta$ and $\eta$ are present under the null. In contrast, such a reparameterization is not required for Example 2 since $\zeta$ is not present under the null.

When the models involved are parametric, a special case when $\eta$ does not depend on $\zeta$ under the null, that is, $\zeta$ is only present under the alternative, has been studied extensively by Andrews and Ploberger (1994); Davies (1977, 1987); King and Shively (1993), and others. Davies (1977) showed that the likelihood ratio test is optimal in the sense that as the significance level of the test tends to zero, its power function approaches that of the optimum test when $\zeta$ is given. These optimality results are very weak and do not provide any guidance regarding the performance of the test in practical applications, where the significance level is fixed, for example, at level 0.05 [Andrews (1999)]. Andrews and Ploberger (1994) studied optimal tests for parametric models using the weighted average power criterion originally introduced by Wald (1943) when studying the likelihood ratio test under regularity conditions, where the model is identifiable under the null. Under loss of identifiability, the likelihood ratio test is generally less powerful than the optimal test in Andrews and Ploberger (1994). These optimal tests possess a Bayesian interpretation, where the weight corresponds to a prior on the nonidentifiable parameter, and are asymptotically equivalent to a Bayesian posterior odds ratio.

In this paper, we adapt the weighted average power criterion [Andrews and Ploberger (1994), Wald (1943)] to construct optimal tests in semiparametric models under loss of identifiability. Our main contribution is to extend the results of Andrews and Ploberger (1994) in at least four directions.

First, Andrews and Ploberger (1994) address only parametric models, as is the case for most of the literature on testing problems with nonidentifiability under the null. Our optimality results are available for semiparametric models, where $\eta$ may be infinite dimensional and $\zeta$ may not be estimable at the usual parametric rate under either the null or the alternative. A semiparametric profile likelihood approach is adopted to reduce the infinite-dimensional model to a finite-dimensional uniformly least-favorable submodel; see Murphy and van der Vaart (2000) for a discussion of profile likelihood in regular settings. We note however, that the idea of uniformly least favorable submodels is a new concept in semiparametric settings, which is not discussed in Murphy and van der Vaart (2000). The development of this concept is both nontrivial and critical to establishing an appropriate optimality criterion for semiparametric models under loss of identifiability.

Second, the results of Andrews and Ploberger (1994) are applicable only for tests where a nuisance parameter (namely $\zeta$) is present only under the alternative. This may not be true in our situation, where a nondegenerate reparameterization may be needed to make $\zeta$ vanish under the null. Furthermore, our tests and the optimality results do not depend on the reparameterization.

Third, Andrews and Ploberger (1994) establish that their test is optimal with respect to local alternatives for $\psi$ involving a multivariate normal prior with singular covariance matrix. In our approach, it is only necessary to specify the prior in the direction of $\beta$, the parameter of interest, and no prior is needed on the remaining parameter $\eta$. This enables us to avoid the singular covariance issue in Andrews and Ploberger (1994).

Fourth, we develop a simple and effective Monte Carlo method of inference for the proposed test statistics.

Adopting a profile likelihood approach has several advantages. First, under the identifiable submodel, the MLE for $\eta$ may converge at a slower rate than the usual $\sqrt{n}$ rate, such as the change-point Cox model with current status data. This makes the theoretical justification based on Taylor expansion of the full likelihood fail. Second, even if the MLE of the nonparametric component converges at the $\sqrt{n}$ rate, semiparametric likelihoods may not be suitably "differentiable," in particular, when such a likelihood contains certain empirical terms, as with, for example, the odds-rate model. Third, handling the remainder terms in a Taylor type expansion is challenging, owing to the presence of the infinite dimensional parameters, and a delicate Banach space analysis is required. Employing the profile likelihood enables us to address these issues rigorously.

1.4. *Organization of paper.* The remainder of the paper is organized as follows. In Section 2, we present the generic testing problem and the model and data assumptions. The optimality results are given in Section 3. We verify that the results hold for the examples in Section 4. In Section 5, we describe a simulation study to evaluate the finite-sample behavior of the proposed tests and to compare its efficiency with some alternative tests for the current status example. In Section 6, we discuss some additional examples without identifiability under the null which are not covered in our current settings and which require further extensions. Proofs are given in Section 7.

## 2. The hypothesis tests and assumptions.

2.1. *The optimal tests.* In this subsection we formulate the tests of hypotheses when the parameters are not identifiable under the null. Let $P_\theta$ denote the probability measure, based on observed data $\tilde{X}_n \equiv (X_1, X_2, \ldots, X_n)$, where $\theta \in \Theta$ and the subscript $n$ is the sample size. As mentioned previously, the parameters $\theta \in \Theta_0$ under the null hypothesis are not identifiable. We assume, as in the examples, that $\theta$ can be partitioned as $(\psi, \zeta)$, with $\zeta$ $q$-dimensional and $\psi$ of arbitrary dimension. We further assume that $\psi$ can be partitioned as $(\beta, \eta)$ so that the null hypothesis can be stated in terms of $\beta$, with the nuisance parameter $\eta$ having arbitrary dimension. The likelihood function of the data is given by $l_n(\theta)$ and the profile likelihood for $\beta$ and $\zeta$ is defined as $pl_n(\beta, \zeta) = \sup_\eta l_n(\beta, \eta, \zeta)$. For the semiparametric model $\{P_{(\beta, \eta, \zeta)}\}$ on a sample space $\mathcal{X}$, we assume $\beta \in \mathbb{R}^p$, $\zeta \in \Xi$, a compact subset of $\mathbb{R}^q$ and $\eta \in \mathcal{H}_\eta$, which is a subset of a Banach space.

The hypotheses to be tested are:

$$(3) \qquad\qquad H_0 : \beta = \beta_0 \quad \text{vs.} \quad H_1 : \beta \neq \beta_0.$$

When $\beta = \beta_0$, the null distribution $P_0$ is unique and the likelihood for a single observation under the null is abbreviated as $l^0$. Let $\pi \equiv (\eta, \zeta)$. The null set of $\pi$ is $\Pi_0$ and its cardinality is the same as that of $\Xi$, which is at least two. $\Theta_0 = \{\beta_0\} \times \Pi_0$. For each $\zeta \in \Xi$, $\eta_0(\zeta) \equiv \{t \in \mathcal{H}_\eta : (t, \zeta) \in \Pi_0\}$ is an interior point of $\mathcal{H}_\eta$. Let $\psi_0(\zeta) \equiv (\beta_0, \eta_0(\zeta))$, and $\theta_0(\zeta) \equiv (\psi_0(\zeta), \zeta)$. Thus, $\Theta_0$ can be represented as $\Theta_0 = \{\theta_0(\zeta) : \zeta \in \Xi\}$.

Before introducing the optimal tests, we need some additional notations for the parameter space and the score and information operators in the semiparametric settings. We denote $\dot{l}_\beta \in L_2^0(P_\theta)$ as the derivative of $\log l_1(\theta)$ with respect to $\beta$ and $\ddot{l}_\beta$ is the second derivative of $\log l_1(\theta)$ with respect to $\beta$. $L_2^0(P_\theta)$ refers to the class of square integrable functions under the measure $P_\theta$ with mean 0. The score operator for $\eta$ is defined as $\dot{l}_\eta$, which is a bounded linear map from $\mathcal{H}_\eta$ to $L_2^0(P_\theta)$ with adjoint operator $\dot{l}_\eta^\star : L_2^0(P_\theta) \mapsto \overline{\mathcal{H}_\eta}$, where $\overline{\mathcal{H}_\eta}$ is the closed linear span of $\mathcal{H}_\eta$. The information operator is $\dot{l}_\eta^\star \dot{l}_\eta : L_2^0(P_\theta) \mapsto L_2^0(P_\theta)$. The efficient score for $\beta$ is the ordinary score function $\dot{l}_\beta$ minus its orthogonal projection onto the closed linear span of the score operator $\dot{l}_\eta$. The efficient information for $\beta$ is $\tilde{I}_\beta = \int \tilde{l}_\beta \tilde{l}'_\beta \, dP_\theta$, which is the asymptotic variance of the efficient score function.

We use the notations $\mathbb{P}_n$ and $\mathbb{G}_n$ for the empirical distribution and the empirical process of the observations. That is, for every measurable function $f$ and probability measure $P$,

$$\mathbb{P}_n f = \frac{1}{n} \sum_{i=1}^n f(X_i), \qquad Pf = \int f \, dP, \qquad \mathbb{G}_n f = \frac{1}{\sqrt{n}} \sum_{i=1}^n (f(X_i) - P(f)).$$

We note that although simultaneous estimation of $\beta$ and $\zeta$ fails under the null due to nonidentifiability, estimation results for $\hat{\beta}_n(\zeta)$, the MLE of $\beta$ at a fixed value of $\zeta$, are often valid under the null. This suggests making inference about $\beta$ using $\hat{\beta}_n(\zeta)$. For fixed $\zeta \in \Xi$, the score, Wald and likelihood ratio test statistics for testing $H_0$ against $H_1$ are given by

$$R_n(\zeta) = \mathbb{P}_n \dot{l}_\beta(\hat{\theta}_0(\zeta))' \{\mathbb{P}_n \dot{l}_\beta \dot{l}'_\beta(\hat{\theta}_0(\zeta))\}^{-1} \mathbb{P}_n \dot{l}_\beta(\hat{\theta}_0(\zeta)),$$

$$W_n(\zeta) = (\hat{\beta}_n(\zeta) - \beta_0)' \hat{\tilde{I}}_\beta(\hat{\theta}_n(\zeta))(\hat{\beta}_n(\zeta) - \beta_0) \quad \text{and}$$

$$LR_n(\zeta) = -2\{l_n(\hat{\theta}_0(\zeta)) - l_n(\hat{\theta}_n(\zeta))\},$$

where $\hat{\theta}_n(\zeta) \equiv (\hat{\beta}_n(\zeta), \hat{\eta}_n(\zeta), \zeta)$ is the unrestricted MLE of $\theta$ at a fixed value of $\zeta$ and $\hat{\theta}_0(\zeta) \equiv (\beta_0, \hat{\eta}_0(\zeta), \zeta)$ is the restricted MLE of $\theta$ for a fixed value of $\zeta$ under the null. $\mathbb{P}_n \dot{l}_\beta(\hat{\theta}_0(\zeta)) = \mathbb{P}_n \dot{l}_\beta(\beta_0, \hat{\eta}_0(\zeta), \zeta)$ is the empirical score function of $\beta$ evaluated at the restricted MLE $\hat{\theta}_0(\zeta)$. $\mathbb{P}_n \dot{l}_\beta(\hat{\theta}_n(\zeta)) = \mathbb{P}_n \dot{l}_\beta(\hat{\beta}_n(\zeta), \hat{\eta}_n(\zeta), \zeta)$ is the empirical score function of $\beta$ evaluated at the unrestricted MLE $\hat{\theta}_n(\zeta)$. The inverse matrix of $\hat{\tilde{I}}_\beta(\hat{\theta}_n(\zeta))$, a consistent estimator of the efficient information $\tilde{I}_\beta$ under the null, estimates the covariance matrix of $\hat{\beta}_n(\zeta)$.

The optimal tests we propose take the form

$$ER_n = (1+c)^{-p/2} \int \exp\left(\frac{1}{2}\frac{c}{1+c} R_n(\zeta)\right) dJ(\zeta),$$

$$EW_n = (1+c)^{-p/2} \int \exp\left(\frac{1}{2}\frac{c}{1+c} W_n(\zeta)\right) dJ(\zeta) \quad \text{and}$$

$$ELR_n = (1+c)^{-p/2} \int \exp\left(\frac{1}{2}\frac{c}{1+c} LR_n(\zeta)\right) dJ(\zeta),$$

where $c > 0$ is a known constant and $J(\cdot)$ is a pre-selected integrable prior on $\zeta$. Their optimality will be discussed in Section 3. We note that, in semiparametric settings, the computation of the efficient information may involve high dimensional maximization and nonparametric smoothing. Then the tests $ER_n$ and $EW_n$ may be computationally harder than $ELR_n$. Hence the likelihood ratio based test $ELR_n$ is more attractive in these settings.

In construction of the optimal tests, understanding and computing $\hat{\theta}_0(\zeta)$, may be complicated due to the dependence of the parameter $\theta_0(\zeta)$ on $\zeta$.

Assuming the existence of the following full rank reparameterization, we can eliminate the dependence between $\eta$ and $\zeta$, thereby easing both the theoretical developments and the computations for the proposed tests.

2.2. *Full rank reparameterization: Breaking the dependence between $\eta$ and $\zeta$.* We assume there exists a map $\phi_\zeta : \mathcal{H}_\eta \mapsto \mathcal{H}_\eta$, which is one-to-one and uniformly Hadamard-differentiable at $\eta$ tangentially to $\mathcal{H}_\eta$ over $\zeta \in \Xi$, that is,

$$\sup_{(\eta + t_n h_n(\zeta), \zeta) \in \Pi_0} \left\| \frac{\phi_\zeta(\eta + t_n h_n(\zeta)) - \phi_\zeta(\eta)}{t_n} - \dot{\phi}_\zeta(\eta)(h(\zeta)) \right\| \to 0,$$

as $\sup_{\zeta \in \Xi} \| h_n(\zeta) - h(\zeta) \| \to 0$, and $t_n \to 0$, where $h(\zeta)$ is in the tangent space of $\mathcal{H}_\eta$ for all $\zeta \in \Xi$ and $\| \cdot \|$ denotes the norm of $\mathcal{H}_\eta$. Its derivative $\dot{\phi}_\zeta$ is one-to-one and continuously invertible uniformly over $\zeta \in \Xi$. That is, there exists a positive constant $c$ such that $\| \dot{\phi}_\zeta(\eta_1(\zeta) - \eta_2(\zeta)) \| \geq c \| \eta_1(\zeta) - \eta_2(\zeta) \|$ for every $\eta_1(\zeta)$ and $\eta_2(\zeta)$ in $\mathcal{H}_\eta$ for all $\zeta \in \Xi$. Let $\overline{\eta} \equiv \phi_\zeta(\eta)$, and $l_1(\beta_0, \overline{\eta}, \zeta)(x) \equiv l_1(\beta_0, \phi_\zeta^{-1}(\overline{\eta}), \zeta)(x) = l^0(x)$, where $\zeta$ vanishes under the null, for all $x$ in $\mathcal{X}$.

This reparameterization does not change the likelihood, that is, the equality $l_1(\beta, \eta(\zeta), \zeta)(x) = l_1(\beta, \overline{\eta}, \zeta)(x)$ holds both under the null and the alternative. Under the null, the likelihood $l_1(\beta_0, \eta_0(\zeta), \zeta) = l_1(\beta_0, \overline{\eta}_0, \zeta)$ for a specific $\overline{\eta}_0$, which does not depend on $\zeta$, and $\zeta$ disappears in the null likelihood. We thus reduce the parameter dimension of the null space from $\Pi_0$ to $\mathcal{H}_\eta$. For Example 2 in the Introduction, $\phi$ can be taken to be the identity and thus the reparameterization is not needed. In contrast, a reparameterization is needed for Example 1. We will give the details later in Section 4.

The reason we assume the existence of such a full rank reparameterization is to eliminate the dependence between $\eta$ and $\zeta$. The issue is that the optimality results are with respect to a perturbation of the parameter $\eta$, which is not well defined in the original space, due to the dependence between parameters $\eta$ and $\zeta$. Subsequent assumptions are built on the new parameterization $\overline{\theta} \equiv (\beta, \overline{\eta}, \zeta)$. However, the results still hold for the original parameterization, since the efficient score and efficient information of $\beta$ are invariant under such reparameterization of $\eta$, as given in the following lemma:

LEMMA 1. *Under the full rank reparameterization, $\tilde{l}_\beta(\theta) = \tilde{\ell}_\beta(\overline{\theta})$, where $\tilde{\ell}_\beta(\overline{\theta})$ is the efficient score of $\beta$ under the new reparameterization. The efficient information matrix is also invariant under these reparameterizations.*

REMARK 1. The full rank reparameterization defined above may not be unique. We will show later in the proof of Theorem 2 that the optimal tests proposed in this paper are invariant to the choice of the full rank reparameterization.

Next we discuss how to construct the optimal tests with the new parameterization, where $\zeta$ vanishes, and $\overline{\psi}$ does not depend on $\zeta$ under the null.

### 2.3. Constructing optimal tests under the full rank reparameterization.
Though $\zeta$ disappears in the likelihood under the null hypothesis, the score and information are still processes indexed by $\zeta$. For fixed $\zeta \in \Xi$, the score, Wald and likelihood ratio test statistics for testing $H_0$ against $H_1$ with the new parameterization can be represented as:

$$R_n(\zeta) = \mathbb{P}_n \dot{\ell}_\beta(\hat{\overline{\psi}}_0, \zeta)' \{\mathbb{P}_n \dot{\ell}_\beta \dot{\ell}_\beta'(\hat{\overline{\psi}}_0, \zeta)\}^{-1} \mathbb{P}_n \dot{\ell}_\beta(\hat{\overline{\psi}}_0, \zeta),$$

$$W_n(\zeta) = (\hat{\beta}_n - \beta_0)' \hat{\tilde{\mathcal{I}}}_\beta(\hat{\overline{\psi}}_n, \zeta)(\hat{\beta}_n - \beta_0) \quad \text{and}$$

$$LR_n(\zeta) = -2\{\ell_n(\hat{\overline{\psi}}_0, \zeta) - \ell_n(\hat{\overline{\psi}}_n, \zeta)\},$$

where $\hat{\overline{\psi}}_n$ is the unrestricted MLE of $\overline{\psi}$ and $\hat{\overline{\psi}}_0$ is the restricted MLE of $\overline{\psi}_0$. $\mathbb{P}_n \dot{\ell}_\beta(\beta_0, \hat{\overline{\eta}}_0, \zeta)$ is the empirical score function of $\beta$ evaluated at the restricted MLE $\hat{\overline{\psi}}_0$. $\mathbb{P}_n \dot{\ell}_\beta(\hat{\beta}_n, \hat{\eta}_n, \zeta)$ is the empirical score function of $\beta$ evaluated at the unrestricted MLE $\hat{\overline{\psi}}_n$. The inverse matrix of $\hat{\tilde{\mathcal{I}}}_\beta(\hat{\overline{\psi}}_n, \zeta)$, a consistent estimator of the efficient information of $\beta$ under the null, estimates the covariance matrix of $\hat{\beta}_n$. It is thus obvious that the optimal tests are invariant with respect to the choice of full rank reparameterizations.

To further study the asymptotic distribution and the optimality of the proposed tests, we need the following assumptions, based on the full rank reparameterization. We note that except assumption C, all other assumptions can also be stated with the original parameterization.

### 2.4. The assumptions based on the reparameterization.
To derive asymptotically optimal tests of $H_0$, we consider local alternatives to $H_0$ of the form $\ell_n(\beta_0 + h/\sqrt{n}, \overline{\eta}, \zeta)$ with $\zeta$ and $\overline{\eta}$ unspecified. The optimality criterion will involve a weighted average power criterion, where the averaging is with respect to an integrable prior $Q_\zeta(h)$ on the values of $h$ in $\mathbb{R}^p$ defining local alternatives and an integrable prior $J(\zeta)$ on $\zeta$. Before formally stating the optimality criterion, we give assumptions on the data and the parameter spaces. The first two assumptions postulate the existence of the prior on local alternatives, $Q_\zeta(h)$.

A1 The efficient information function of $\beta$ evaluated at $(\overline{\psi}_0, \zeta)$, $\tilde{\mathcal{I}}_\beta(\overline{\psi}_0, \zeta)$, is uniformly continuous in $\beta$ and $\zeta$ over $B_0 \times \Xi$, where $B_0$ is some neighborhood of $\beta_0$. Furthermore, $\tilde{\mathcal{I}}_\beta(\overline{\psi}_0, \zeta)$ is uniformly positive definite over $\zeta \in \Xi$, that is $\inf_{\zeta \in \Xi} \lambda_{\min}\{\tilde{\mathcal{I}}_\beta(\overline{\psi}_0, \zeta)\} > 0$, where $\lambda_{\min}(C)$ is the smallest eigenvalue of the matrix $C$.

A2 $Q_\zeta$ is a normal measure with mean $\beta_0$ and variance $c\tilde{\mathcal{I}}_\beta^{-1}(\overline{\psi}_0, \zeta)$ for $\zeta \in \Xi$,
   where $c > 0$ is a scalar constant.

Assumptions A1 and A2 are analogous to Assumptions 1(e), 1(f) and 4 of
Andrews and Ploberger (1994), although there are fundamental differences.
Andrews and Ploberger (1994) work directly by building on the full para-
metric likelihood and their assumptions refer to the information matrix for
all parameters. Furthermore, their optimality results are defined in terms
of local alternatives for $\psi$, where the prior is a multivariate normal with
singular covariance matrix. Our assumptions A1 and A2 are only for the pa-
rameter of interest, $\beta$, with no prior assumptions needed for $\eta$ under either
the null or the alternative.

The next set of conditions assumes the existence of a uniformly least-
favorable submodel. This submodel can be viewed as a "uniform" version of
the least favorable submodel discussed in Murphy and van der Vaart (2000):
the convergence rate of the nuisance parameter now is in the "uniform"
sense, and the efficient score and the efficient information possess Donsker
and Glivenko–Cantelli properties with "larger" index sets, respectively. When
the set of $\zeta$, $\Xi$, is a singleton, this new submodel concept reduces to the ordi-
nary least favorable submodel. The development of this concept is critical to
establishing an appropriate optimality criterion for general semiparametric
models under loss of identifiability. Here are the needed assumptions:

B1 There exists a map $t \mapsto f_t$ from a fixed neighborhood of $\beta_0$ into $\mathcal{H}_\eta$, such
   that the map $t \mapsto \ell(t, \overline{\theta})$ defined by $\ell(t, \overline{\theta}) \equiv \ell_1(t, f_t, \zeta)$ is twice contin-
   uously differentiable. Let $\dot{\ell}(t, \overline{\theta})$ and $\ddot{\ell}(t, \overline{\theta})$ denote the derivatives with
   respect to $t$. The submodel with parameters $(t, f_t, \zeta)$ passes through $\overline{\eta}$
   at $t = \beta$, that is, $f_\beta(\beta, \overline{\eta}, \zeta) = \overline{\eta}$ for all $\zeta \in \Xi$.
B2 The submodel is uniformly least-favorable at $\overline{\psi}_0 = (\beta_0, \overline{\eta}_0)$ and $\zeta$ for
   estimating $\beta_0$ in the sense that $\dot{\ell}(\beta_0, \overline{\psi}_0, \zeta) = \tilde{\ell}_\beta(\overline{\psi}_0, \zeta)$. As $(t, \beta, \overline{\eta}) \to$
   $(\beta_0, \beta_0, \overline{\eta}_0)$, we assume that $\sup_{\zeta \in \Xi} \|\dot{\ell}(t, \overline{\psi}, \zeta) - \tilde{\ell}_\beta(\overline{\psi}_0, \zeta)\| = o_{P_0}(1)$ and
   $\sup_{\zeta \in \Xi} \|\ddot{\ell}(t, \overline{\psi}, \zeta) - \ddot{\ell}(\beta_0, \overline{\psi}_0, \zeta)\| = o_{P_0}(1)$. In the sequel, we let $o_{\overline{P}}^\Xi$ denote
   a quantity going to zero in probability, under $P$, uniformly over the set
   $\Xi$.
B3 We assume that $\hat{\overline{\psi}}_0$, the restricted MLE of $\overline{\psi}$ under the null, satisfies
   $\hat{\overline{\psi}}_0 = \overline{\psi}_0 + o_{P_0}(1)$. The unrestricted MLE $\hat{\overline{\psi}}_n(\zeta) = \overline{\psi}_0 + o_{\overline{P}_0}^\Xi(1)$. More-
   over, let $\hat{\overline{\eta}}_\beta(\zeta) \equiv \arg\max_{\overline{\eta}} \ell_n(\beta, \overline{\eta}, \zeta)$, that is, $pl_n(\beta, \zeta) = \ell_n(\beta, \hat{\overline{\eta}}_\beta(\zeta), \zeta)$.
   Assume that for any random sequences $\tilde{\beta}_n \to_{P_0} \beta_0$, we have $\hat{\overline{\eta}}_{\tilde{\beta}_n}(\zeta) =$
   $\overline{\eta}_0 + o_{\overline{P}_0}^\Xi(1)$ and the following uniform "no-bias" condition holds:

$$(4) \qquad P_0 \dot{\ell}(\beta_0, \tilde{\beta}_n, \hat{\overline{\eta}}_{\tilde{\beta}_n}(\zeta), \zeta) = o_{\overline{P}_0}^\Xi(\|\tilde{\beta}_n - \beta_0\| + n^{-1/2}).$$

B4 There exist neighborhoods $U$ of $\beta_0$ and $V$ of $\overline{\psi}_0$, such that the class of functions $\{\dot{\ell}(t, \overline{\psi}, \zeta) : t \in U, \overline{\psi} \in V, \zeta \in \Xi\}$ is $P_0$-Donsker with square integrable envelope function and the class of functions $\{\ddot{\ell}(t, \overline{\psi}, \zeta) : t \in U, \overline{\psi} \in V, \zeta \in \Xi\}$ is $P_0$-Glivenko–Cantelli and is bounded in $L_1(P_0)$, where $L_1(P_\theta)$ refers to the class of integrable functions under $P_\theta$.

Assumptions B1–B4 set the stage for the quadratic expansion of the profile likelihood and the derivation of the optimality properties of the proposed tests. Note that these assumptions can also be built on the original parameterization, but we use the new parameterization for ease of presentation. Since our formulation includes parametric models as special cases, the existence of a uniformly least-favorable submodel in our set-up covers all situations considered by Andrews and Ploberger (1994).

Compared with Andrews and Ploberger (1994), we have a stronger form of the unbiasedness condition and stronger requirements on the consistency of the estimators for the expansion of the profile likelihood. This is partly due to the more general structure of the semiparametric model. As in assumption B3, we require that if $\tilde{\beta}_n$ is any sequence of estimators consistent for $\beta_0$, $\hat{\overline{\eta}}_\zeta(\tilde{\beta}_n)$ must be consistent for $\overline{\eta}_0$, the true value of the nuisance parameter $\overline{\eta}$, uniformly over $\Xi$. In Andrews and Ploberger (1994), consistency is only needed for the unconstrained MLE (assumption 2) and the constrained MLE under the null hypothesis (assumption 3).

To evaluate the local asymptotic distribution of the proposed tests, we require differentiability in quadratic mean (DQM) of the parameters $\overline{\psi}$, as stated in the following assumption C, which is commonly used to evaluate the local power. It will be verified for the two examples presented in the introduction. Unlike assumptions B1–B4, the full rank reparameterization is indispensable in assumption C:

C Differentiability in quadratic mean of the parameter $\overline{\psi}$. A perturbation of $\overline{\psi}$ in its domain is $\overline{\psi}_t = \overline{\psi}_0 + th + o(1)$, where $h \equiv (h_\beta, h_{\overline{\eta}})$, $h_\beta \in \mathbb{R}^p$ and $h_{\overline{\eta}} \in \overline{\mathcal{H}}_\eta$. The DQM condition for $\overline{\psi}_0$ with respect to the collection of paths $\{\overline{\psi}_t\}$ is:

$$\int \left[ \frac{(dP_{\overline{\psi}_t, \zeta})^{1/2} - (dP_0)^{1/2}}{t} - \frac{1}{2}(A_\zeta h) \, dP_0^{1/2} \right]^2 \to 0, \qquad \text{as } t \to 0,$$

for all $\zeta \in \Xi$, where $A_\zeta$ is a bounded linear operator defined on $\mathbb{R}^p \times \overline{\mathcal{H}}_\eta$ and takes values in $L_2^0(P_\theta)$.

Differentiability in quadratic mean implies that the range of $A_\zeta$ is contained in $L_2^0(P_\theta)$. Note that $A_\zeta h = (\partial/\partial t)\ell_1(\overline{\psi}_t, \zeta)|_{t=0}$, following similar arguments as in Kosorok and Song (2007), where $h = (h_\beta, h_{\overline{\eta}})$. We define $A_\zeta$ to be given by $A_\zeta(h_\beta, h_{\overline{\eta}}) = \dot{\ell}'_\beta(\overline{\psi}, \zeta)h_\beta + \dot{\ell}_{\overline{\eta}}(\overline{\psi}, \zeta)h_{\overline{\eta}}$, where $\dot{\ell}_\beta$ and $\dot{\ell}_{\overline{\eta}}$ are the score operators for $\beta$ and $\overline{\eta}$, respectively. Moreover, $\mathbb{R}^p \times \overline{\mathcal{H}}_\eta$ is a Hilbert space with $\|\cdot\|$ denoting its norm and $\langle \cdot, \cdot \rangle$ denoting its inner product. Since in parametric settings, twice continuous differentiability implies DQM [Pollard (1995)], this assumption is weaker than Assumption 1(c) in Andrews and Ploberger (1994).

**3. Main results.** This section includes several main results. The first one gives the asymptotic null distribution of the proposed tests.

3.1. *The distributions of the test statistics under the null.* To establish the asymptotic null distribution of the test statistics, a key result about the uniform profile likelihood expansion is summarized in the following lemma.

LEMMA 2. *Under assumptions* A–C, *for any random sequence* $\tilde{\beta}_n \to_{P_0} \beta_0$,

$$
\begin{aligned}
\log pl_n(\tilde{\beta}_n, \zeta) = {} & \log pl_n(\beta_0, \zeta) + n(\tilde{\beta}_n - \beta_0)' \mathbb{P}_n \tilde{l}_\beta(\theta_0(\zeta)) \\
& - \tfrac{1}{2} n(\tilde{\beta}_n - \beta_0)' \tilde{I}_\beta(\theta_0(\zeta))(\tilde{\beta}_n - \beta_0) \\
& + o_{P_0}^{\overline{\overline{=}}}(\sqrt{n}\|\tilde{\beta}_n - \beta_0\| + 1)^2.
\end{aligned}
\tag{5}
$$

Lemma 2 enables us to establish the asymptotic equivalence of these test statistics and their asymptotic distributions:

THEOREM 1. *Under assumptions* A–C, $ELR_n = EW_n + o_{P_0}(1) = ER_n + o_{P_0}(1) \to_d e\chi(c)$, *where*

$$
e\chi(c) = (1+c)^{-p/2} \int \exp\left( \frac{1}{2} \frac{c}{1+c} \mathbb{G}'(\theta_0(\zeta)) \tilde{I}_\beta^{-1}(\theta_0(\zeta)) \mathbb{G}(\theta_0(\zeta)) \right) dJ(\zeta),
$$

*and* $\mathbb{G}(\theta_0(\zeta))$ *is the limiting process of* $\mathbb{G}_n \tilde{l}_\beta(\theta_0(\zeta))$, *which is a mean zero Gaussian process with variance function* $\sigma^2(\zeta) = \tilde{I}_\beta(\theta_0(\zeta))$ *indexed by* $\zeta$ *and with covariance function* $\sigma^2(\zeta_1, \zeta_2) = P_0\{\tilde{l}_\beta(\theta_0(\zeta_1))\tilde{l}_\beta(\theta_0(\zeta_2))'\}$, *indexed by* $\zeta_1$ *and* $\zeta_2$, *for* $\zeta$, $\zeta_1$ *and* $\zeta_2 \in \Xi$.

REMARK 2. When $J(\cdot)$ does not correspond to a prior on $\zeta$, corresponding rather to a weight function lacking a probabilistic interpretation, then the results in Theorem 1 will generally hold, although the test may no longer possess the optimality discussed in the sequel. Theorem 1 should also hold if $Q_\zeta(h)$ is not a prior distribution, corresponding rather to a weight function on local alternatives for $\beta$. This robustness indicates that the tests are generally valid under loss of identifiability, yielding a large class of test statistics, with the optimal test being a member of this class.

We note that Theorem 1 only holds for normal weight $Q_\zeta^c$, which corresponds to the uniform least favorable direction. As indicated in the proof of Theorem 1, the normal weight function $Q_\zeta(\cdot)$ is integrated out, hence does not appear in the test with the original form. Subsequently, the optimal tests depend on the weight function $Q_\zeta(\cdot)$ only through the scalar $c$. The larger $c$ is, the more weight is given to alternatives for which $\beta$ is large. For example, for a test of the change-point model, larger values of $c$ correspond to greater weight being given to larger changes. In the special case where $J(\zeta)$ is a pointmass at a single value $\zeta_0$, the optimal test rejects if and only if $LR(\zeta_0)$ exceeds some constant (i.e., the optimal test equals the standard score test for fixed $\zeta_0$) and the optimal test is independent of $c$. When $J(\zeta)$ is not a pointmass distribution, however, the optimal test $ELR_n$ depends on $c$. The larger $c$ is, the more power is directed at alternatives for which $\beta$ is large.

The limit as $c \to 0$ of the $2(ELR_n - 1)/c$ statistic is equal to the "average score" statistic $\int LR_n(\zeta) \, dJ(\zeta)$, which is the limit of the $ELR$ statistics that are designed for alternatives that are very close to the null hypothesis. At the other extreme, the limit as $c \to \infty$ is $\log \int \exp(LR_n(\zeta)/2) \, dJ(\zeta)$. Thus for testing against more distant alternatives, the optimal test statistic is still of an average exponential form.

If the constant $c/(1 + c)$ which appears in the definition of $ELR_n$ is replaced by a constant $r > 0$, then the limit as $r \to \infty$ of $ELR_n$ is the likelihood ratio test, equivalently, the "sup score" statistic studied in Kosorok and Song (2007). Hence, the sup score test is designed for distant alternatives, but is of a more extreme form than the optimal exponential test, since the latter requires $r < 1$. It can be easily shown as a corollary to Theorem 1 that the usual likelihood ratio, Wald and score tests have the following distribution:

COROLLARY 1. *Under the null hypotheses and assumptions* A–C, $\sup_\zeta LR_n(\zeta) = \sup_\zeta W_n(\zeta) + o_{P_0}(1) = \sup_\zeta R_n(\zeta) + o_{P_0}(1) \to_d \chi$, *with* $\chi = \sup_\zeta \mathbb{G}'(\theta_0(\zeta))\tilde{I}_\beta^{-1}(\theta_0(\zeta))\mathbb{G}(\theta_0(\zeta))$.

3.2. *Optimality of the proposed tests.* The second main result of this paper is the optimality property of the proposed tests. Following assumptions in Section 2, we consider local alternatives $\beta = \beta_0 + h_\beta/\sqrt{n} + o(n^{-1/2})$ for $h_\beta \in \mathbb{R}^p$ with prior distribution $Q_\zeta(h_\beta)$ on the local alternative direction $h_\beta$ and prior distribution $J(\zeta)$ on the nonidentifiable parameter $\zeta$. The optimality result is as follows:

THEOREM 2. *Under assumptions* A–C, *the test statistics in Theorem 1 are asymptotically uniformly most powerful for testing* $H_0 : \beta = \beta_0$ *against the contiguous alternative*

$$\int dP_{\bar\psi_0 + h/\sqrt{n} + o(n^{-1/2}), \zeta}^n \, dQ_\zeta(h_\beta) \, dJ(\zeta),$$

*where* $h \equiv (h_\beta, h_\eta(\zeta))$, $h_\eta(\zeta) \equiv \tilde{q}'_\zeta h_\beta$ *and where* $\tilde{q}_\zeta \equiv -(\dot{\ell}^\star_\eta \dot{\ell}_\eta)^- \dot{\ell}^\star_\eta \dot{\ell}_\beta(\overline{\psi}, \zeta)$ *is the uniformly least-favorable direction indexed by* $\zeta$. *Moreover, this optimality result is invariant under the choice of reparameterization.*

Theorem 2 also implies that the proposed tests have the greatest weighted average power asymptotically in the class of all tests of asymptotic significance level $\alpha$, against the alternative $P^n_{\overline{\psi}_0 + h/\sqrt{n} + o(n^{-1/2}), \zeta}$. That is, they maximize

$$\overline{\lim_{n \to \infty}} \int P(\phi_n \text{ rejects} | \overline{\psi}_0 + h/\sqrt{n} + o(n^{-1/2}), \zeta) \, dQ_\zeta(h_\beta) \, dJ(\zeta)$$

over all tests $\phi_n$ of asymptotic level $\alpha$.

Our optimality results are under alternatives $\beta_0 + h_\beta/\sqrt{n} + o(n^{-1/2})$, with nonsingular normal weights on $h_\beta$. Our weights on $h_\beta$ are precisely Andrews and Ploberger's [2] weights projected onto the parameter space that is of interest. Thus, our results and Andrews and Ploberger's are consistent.

We now discuss the choice of the direction $q_\zeta$, the priors $Q_\zeta(\cdot)$ and $J(\cdot)$. By the Neyman–Pearson lemma, for any appropriate prior distributions $Q_\zeta(\cdot)$ and $J(\cdot)$ and any known directions $q_\zeta$, a UMP test for testing $H_0 : \beta = \beta_0$ against the contiguous alternative $\int dP^n_{\overline{\psi}_0 + h/\sqrt{n} + o(n^{-1/2}), \zeta} \, dQ_\zeta(h_\beta) \, dJ(\zeta)$, where $h \equiv (h_\beta, h_\eta(\zeta))$, $h_\eta(\zeta) = q'_\zeta h_\beta$ is defined by

$$\gamma_n = \begin{cases} 1, & \text{if } QLR_n > k_{\alpha n}, \\ \lambda_n, & \text{if } QLR_n = k_{\alpha n}, \\ 0, & \text{if } QLR_n < k_{\alpha n}, \end{cases}$$

where $k_{\alpha n} > 0$, $\lambda_n \in [0,1]$ are constants such that the rejection probability is $\alpha$ under the null and

$$QLR_n = \frac{\int l_n(\overline{\psi}_0 + h/\sqrt{n} + o(n^{-1/2}), \zeta) \, dQ_\zeta(h_\beta) \, dJ(\zeta)}{l^0_n}.$$

We have the following result:

COROLLARY 2. *Under assumptions* A–C, *the null hypothesis and the contiguous alternatives,*

$$QLR_n = (1+c)^{-p/2} \int \exp\left(\frac{1}{2} \frac{c}{1+c} LR_n(\zeta)\right) W(q_\zeta, \zeta) \, dJ(\zeta) + o_p(1),$$

*where* $W(q_\zeta, \zeta) \leq 1$ *is defined in equation (17) in Section 7 below. When* $q_\zeta = \tilde{q}_\zeta$, $W(\tilde{q}_\zeta, \zeta) = 1$ *and* $QLR_n = ELR_n + o_{P_0}(1)$.

As the alternatives we consider are contiguous to the null, in each direction $q_\zeta$, which indexes $QLR_n$, there exists a consistent estimator $\tilde{\eta}_n(q_\zeta)$ of $\overline{\eta}_0$ by the convolution theorem, provided certain conditions hold. The optimal tests can thus be built on $\bar{\eta}_n(q_\zeta)$.

In applications with composite hypotheses, where $q_\zeta$ is unknown, there may not exist a direction which can maximize the power over all directions [Bickel, Ritov and Stoker (2006)]. In a regular testing problem, where all parameters are identifiable, it can be shown that the likelihood ratio test, which is built on the uniformly least-favorable direction, will maximize the minimum power of all directions of the alternatives, over all the test based directions. In our nonregular testing problem, the situation is further complicated, since the power depends on the covariance structure of $\mathbb{G}(\theta_0(\zeta))$. It is not clear if the maximin property still holds in our problem. We note that, however, our tests can be interpreted as the "maximum direction" test. Moreover, since the power of the test is not affected by multiplying by a constant in $QLR_n$, we can standardize $W(q_\zeta, \zeta)\,dJ(\zeta)$ to obtain $d\tilde{J}(\zeta)$, which is a probability measure on $\zeta$. Then the question of the optimal choice of both $q_\zeta$ and $J(\zeta)$ reduces to the question of the optimal choice of $\tilde{J}(\zeta)$. Hence, without loss of generality we can replace $q_\zeta$ with $\tilde{q}_\zeta$. For this reason, we should choose $q_\zeta = \tilde{q}_\zeta$ and focus on the choice of $Q_\zeta(\cdot)$ and $\tilde{J}(\cdot)$ for optimization.

One reason we use the normal weight for $Q_\zeta$ in this paper is to facilitate a comparison with Andrews and Ploberger (1994). Using the normal prior with covariance matrix proportional to the efficient information matrix also leads to a significant simplification of the representation of the test statistics, since many terms cancel in the proof of Theorem 1. However we note that the choice of $Q_\zeta(\cdot)$ is not limited to the normal weight studied in this paper, as indicated in the proof of Theorem 2. More general choices of the priors $Q_\zeta(\cdot)$ and $J(\cdot)$ merit future consideration, but this is beyond the scope of the current paper.

The optimality of the likelihood ratio statistics with loss of identifiability under the null for semiparametric models is of potential interest. Similar to the likelihood ratio test under loss of identifiability with parametric models [Andrews and Ploberger (1994)], in the semiparametric setting, the profile likelihood ratio statistic is not of the optimal average exponential form. It can be shown to be a limit of an average exponential test, but only if a parameter is pushed beyond an admissible boundary, as noted by Andrews and Ploberger (1995) in the parametric case.

3.3. *The distributions of the test statistics under local alternatives.* To gain insight into the power of the optimal tests in practice, it is worthwhile to study their asymptotic distributions under local alternatives. In the following two theorems, Theorem 3 gives the asymptotic distribution for fixed local alternatives $P^n_{\psi_0 + h/\sqrt{n} + o(n^{-1/2}), \zeta_1}$, while Theorem 4 gives the asymptotic distribution for random local alternatives $\int dP^n_{\psi_0 + h/\sqrt{n} + o(n^{-1/2}), \zeta}\,dQ_\zeta(h_\beta)\,dJ(\zeta)$.

As shown in the theorems, the distributions depend on the form of the alternative, which will depend in part on the specifics of the application. These results also usually depend on the prior distributions $J(\cdot)$ and $Q_\zeta(\cdot)$, for both fixed alternatives and random alternatives, although in different manners.

THEOREM 3. *Under local alternatives $P^n_{\psi_0 + h/\sqrt{n} + o(n^{-1/2}), \zeta_1}$ and assumptions A–C, $ELR_n = EW_n + o_p(1) = ER_n + o_p(1) \to_d f\chi(c)$, with*

$$f\chi(c) = (1+c)^{-p/2} \int \exp\left[\frac{1}{2}\frac{c}{1+c}\{\mathbb{G}(\theta_0(\zeta)) + \nu_\star(h_\beta, \zeta, \zeta_1)\}'\right.$$
$$\left. \times \tilde{I}_\beta^{-1}(\theta_0(\zeta))\{\mathbb{G}(\theta_0(\zeta)) + \nu_\star(h_\beta, \zeta, \zeta_1)\}\right] dJ(\zeta),$$

*where $\nu_\star(h_\beta, \zeta, \zeta_1) \equiv P_0 \tilde{l}_\beta(\theta_0(\zeta))\tilde{l}_\beta(\theta_0(\zeta_1))' h_\beta$.*

Now we establish the asymptotic distribution of the test statistics under the alternative $\int dP^n_{\psi_0 + h/\sqrt{n} + o(n^{-1/2}), \zeta} dQ_\zeta(h_\beta) dJ(\zeta)$.

THEOREM 4. *Under assumptions A–C and the local alternative $\int dP^n_{\psi_0 + h/\sqrt{n} + o(n^{-1/2}), \zeta} dQ_\zeta(h_\beta) dJ(\zeta)$, $ELR_n = EW_n + o_p(1) = ER_n + o_p(1) \to_d r\chi(c)$, where $r\chi(c)$ is a real random variable such that its cumulative distribution function $\Pr(r\chi(c) \leq t) = P_0[1\{e\chi(c) \leq t\}e\chi(c)]$.*

3.4. *Monte Carlo computation and inference.* Although we have obtained the asymptotic distributions of the test statistics, these distributions generally have complicated analytic forms which depend on the values of unknown nuisance parameters. We now introduce a weighted bootstrap method to obtain the asymptotically valid critical values of $e\chi(c)$. This method does not require explicit evaluation of the limiting distribution, thereby avoiding the numerical difficulties inherent in such an evaluation.

We first generate $n$ i.i.d. positive random variables $\kappa_1, \ldots, \kappa_n$, with mean $0 < \mu_\kappa < \infty$, variance $0 < \sigma_\kappa^2 < \infty$ and with $\int_0^\infty \sqrt{P(\kappa_1 > u)}\, du < \infty$. Next, we divide each weight by the sample average of the weights $\bar{\kappa}$, to obtain "standardized weights" $\kappa_1^\circ, \ldots, \kappa_n^\circ$ which sum to $n$. For a real, measurable function $f$, define the weighted empirical measure $\mathbb{P}_n^\circ f \equiv n^{-1} \sum_{i=1}^n \kappa_i^\circ f(X_i)$. Let $\hat{\psi}_n^\circ(\zeta) = (\hat{\beta}_n^\circ(\zeta), \hat{\eta}_n^\circ(\zeta))$ denote the maximizer of $l_n^\circ(\psi, \zeta)$ over $\psi \in \Psi$ at fixed $\zeta \in \Xi$, where $l_n^\circ$ is obtained by replacing $\mathbb{P}_n$ with $\mathbb{P}_n^\circ$ in the definition of $l_n$. Similarly, let $\hat{\psi}_0^\circ(\zeta) = (\hat{\beta}_0^\circ(\zeta), \hat{\eta}_0^\circ(\zeta))$ denote the maximizer of $(l_n^0)^\circ(\psi, \zeta)$ over $\psi \in \Psi$ at fixed $\zeta \in \Xi$, where $(l_n^0)^\circ$ is obtained by replacing $\mathbb{P}_n$ with $\mathbb{P}_n^\circ$ in the definition of $l_n^0$, the log likelihood under the null. Now repeat the bootstrap procedure a large number of times $\tilde{M}_n$ and compute the differences of the bootstrapped unrestricted MLE and restricted

MLE of $\beta$: $d\hat{\beta}_k^\circ(\zeta) = \hat{\beta}_{n,k}^\circ(\zeta) - \hat{\beta}_{0,k}^\circ(\zeta)$, $k = 1, \ldots, M_n$, as processes of $\zeta$. Note that we are allowing the number of bootstraps to depend on $n$. Define $\zeta \mapsto \hat{\mu}_n(\zeta) \equiv \tilde{M}_n^{-1} \sum_{k=1}^{\tilde{M}_n} d\hat{\beta}_k^\circ(\zeta)$ and let

$$\zeta \mapsto \hat{V}_n(\zeta) = \tilde{M}_n^{-1} \sum_{k=1}^{\tilde{M}_n} (d\hat{\beta}_{1,k}^\circ(\zeta) - \hat{\mu}_n(\zeta))(d\hat{\beta}_{1,k}^\circ(\zeta) - \hat{\mu}_n(\zeta))'.$$

To estimate critical values, we compute the standardized bootstrap test statistics

$$\begin{aligned} T_{n,k}^\circ \equiv (1+c)^{-p/2} \int \exp\Bigg[ \frac{1}{2} \frac{c}{1+c} \{ (d\hat{\beta}_{1,k}^\circ(\zeta) - \hat{\mu}_n(\zeta))' \\ \times \hat{V}_n^{-1}(\zeta)(d\hat{\beta}_{1,k}^\circ(\zeta) - \hat{\mu}_n(\zeta)) \} \Bigg] \, dJ(\zeta), \end{aligned}$$

for $1 \le k \le \tilde{M}_n$. For a test of size $\alpha$, we compare the observed test statistics with the $(1-\alpha)$th quantile of the corresponding $\tilde{M}_n$ standardized bootstrap statistics. The reason we subtract off the mean is to ensure that we obtain a valid approximation to the null distribution when the null hypothesis may not be true. If not, then there may be loss of power, although the type I error rate will still be controlled when the null is true. The proof of the bootstrap validity can be built upon the proof of Theorems 7 and 8 in Kosorok and Song (2007). We omit the details.

**4. Examples.** In this section, we study the two examples in the introduction to illustrate the two types of nonidentifiability settings, one where a nuisance parameter is present under the null and one where it is not. These examples demonstrate important differences in how the full rank reparameterizations and uniformly least favorable submodels are defined in the two settings. We present Example 2 first because a reparameterization is not required, simplifying the presentation.

4.1. *Example 2 revisited: Change-point regression for current status data.* In the change-point Cox model with current status data, a test of the existence of a threshold effect corresponds to a test of the null $H_0 : \beta = 0$. The change-point parameter $\zeta$ is present only under the alternative. Hence it suffices to take $\phi_\zeta$ as the identity map.

We make the following assumptions and will argue that the assumptions in Section 2 can be checked under these assumptions. Given $Z$, $T$ and $V$ are independent, $Z$ belongs to a compact subset of $\mathbb{R}$. The change-point parameter $\zeta \in [a, b]$, for some known $-\infty < a < b < \infty$ with $\Pr(Z < a) > 0$ and $\Pr(Z > b) > 0$. Assume $P(\text{Var}(Z|V)) > 0$, which guarantees that, as we will show later, the efficient information $\tilde{I}_\beta(\theta_0(\zeta))$ is positive definite

uniformly over $\zeta \in \Xi$. The Lebesgue density of $V$ is positive and continuous on its support $[\sigma, \tau]$ with $0 < \sigma < \tau < \infty$. The baseline hazard function $\Lambda$ is continuously differentiable at $[\sigma, \tau]$, with derivative that is bounded away from 0 and satisfies $\Lambda_0(\sigma) > 0$, $\Lambda_0(\tau) < M$, for some known $M$. We let $\mathcal{H}_\Lambda$ denote a set of nondecreasing cadlag functions $\Lambda$ on $[\sigma, \tau]$ with $\Lambda(\tau) \le M$.

The likelihood function equals (2) with $f_{V,Z}(v, z)$ removed, because it can be absorbed into the underlying measure on the sample space. The log-likelihood for a single observation $\log l_1(\theta)$ takes the form $\log l_1(\theta) = \delta \log[1 - \exp\{-\Lambda(v) \times \exp(r_\gamma(z))\}] - (1 - \delta) \exp(r_\gamma(z)) \Lambda(v)$. Define $Z(\zeta) \equiv (1\{Z > \zeta\}, Z1\{Z > \zeta\}, Z)$ and note that with such a data representation we can adopt much material in the literature and hence simplify our arguments.

To define a uniformly least-favorable submodel in $\beta$, we take two steps. For Step 1, we calculate scores for $\xi$ and $\Lambda$. The score function for $\xi$ is $\dot{l}_\xi(x) = z(\zeta)\Lambda(v)Q(x; \theta)$ with

$$Q(x; \theta) = e^{r_\gamma(z)} \left[ \delta \frac{e^{-e^{r_\gamma(z)}\Lambda(v)}}{1 - e^{-e^{r_\gamma(z)}\Lambda(v)}} - (1 - \delta) \right].$$

The score operator for $\Lambda$ along $\Lambda_t = \Lambda + th$ with $t \ge 0$ and $h$ a nondecreasing nonnegative right continuous function, is given by

$$\dot{l}_\Lambda(h)(x) = \frac{\partial}{\partial t} \log p(x; \gamma, \Lambda_t)|_{t=0} = h(v)Q(x; \theta).$$

We project $\dot{l}_\xi(X)$ onto the space generated by $\dot{l}_\Lambda$. That is, we need to find a function $h_\zeta^\star(V) \in H_\Lambda$ such that $\dot{l}_\xi - \dot{l}_\Lambda(h_\zeta^\star) \perp \dot{l}_\Lambda(h)$, for all $h \in \mathcal{H}_\Lambda$, which is equivalent to solving the least squares problem $P_\theta \|\dot{l}_\xi - \dot{l}_\Lambda h\|^2$. The solution under the null is $h_\zeta^\star(V) \equiv \Lambda_0(V)h_\zeta^{\star\star}(V)$, where $h_\zeta^{\star\star} = P(Z(\zeta)Q^2(X; \psi))/P(Q^2(X; \theta))$, which is assumed to possess a version that is differentiable componentwise with the derivatives being bounded on $[\sigma, \tau]$ uniformly over $\zeta \in \Xi$. It can be shown that $\Lambda_t(\theta)$ is indeed a hazard function when $t$ is sufficiently close to $\xi$.

The uniformly least-favorable direction for $\xi$ is $\Lambda_t(\theta) = \Lambda + (\xi - t)'\varphi(\Lambda)h_\zeta^{\star\star} \circ \Lambda_0^{-1} \circ \Lambda$. Here $\varphi$ is a function mapping $[0, M]$ into $[0, \infty)$ such that $\varphi(y) = y$ on $[\Lambda_0(\sigma), \Lambda_0(\tau)]$ and the function $y \mapsto \varphi(y)/y$ is Lipschitz and $\varphi(y) \le c(\min(y, M - y))$ for a sufficiently large constant $c$. The efficient score for $\xi$ for this uniformly least-favorable submodel is given by:

$$\tilde{l}_\xi(x; t, \theta) = \left[ z - \frac{\varphi(\Lambda)(v)}{\Lambda_t(\theta)(v)} h_\zeta^{\star\star} \circ \Lambda_0^{-1} \circ \Lambda(v) \right] \Lambda_t(\theta)(v) Q(x; t, \Lambda_t(\theta)).$$

$\Lambda_0^{-1}$ may be extended to $[0, \infty)$ by setting $\Lambda_0^{-1}(u) = \sigma$ for $u \le \Lambda_0(\sigma)$ and $\Lambda_0^{-1}(u) = \tau$ for $u > \Lambda_0(\tau)$.

For Step 2, we next project $\dot{l}_\beta(x)$ onto the space generated by $\tilde{l}_\xi$. The efficient score function for $\beta$, $\tilde{l}_\beta$, is the first two coordinates of $\tilde{l}_\xi$ minus its projection on the remaining coordinates of $\tilde{l}_\xi$. Since $\tilde{l}_\xi$ lies in a finite-dimensional space, the projection path has a matrix representation. The efficient information for $\xi$, $\tilde{I}_\xi$, can be partitioned as a two-by-two block matrix, with $\tilde{I}_\xi^{11}(\theta)$ denoting its first two-by-two principle submatrix, and so on. We define $\nu'_\theta = (1, -(\tilde{I}_\xi^{22})^{-1}\tilde{I}_\xi^{21})$, and $\xi_t(\theta) = \xi - (\beta - t)\nu_\theta$. We also define $\Lambda_t(\theta) = \Lambda + (\xi_t(\theta) - t)'\varphi(\Lambda)h_\xi^{\star\star} \circ \Lambda_0^{-1} \circ \Lambda$.

Now we use the uniform least-favorable path $t \mapsto (\xi_t(\theta), \Lambda_t(\theta))$ in the parameter space for the nuisance parameter $\eta \equiv (\alpha, \Lambda)$. This leads to $l(t, \beta, \alpha, \Lambda) = \log l(\xi_t(\theta), \Lambda_t(\theta))$. This submodel is least favorable at $(\xi_0, \Lambda_0)$ uniformly over $\zeta \in \Xi$ since $\partial/\partial t|_{t=\beta_0} l(t, \beta_0, \alpha, \Lambda) = \nu'_\theta \tilde{l}_\xi$, whereas $\nu'_\theta \tilde{l}_\xi = \tilde{l}_\beta$. The efficient information matrix for $\beta$ is, $\tilde{I}_\beta = \tilde{I}_\xi^{11} - \tilde{I}_\xi^{11}(\tilde{I}_\xi^{22})^{-1}\tilde{I}_\xi^{21}(\theta)$. The remainder of assumption B4 can be verified by standard empirical process arguments.

To verify assumption A1 in Section 2, it suffices to show that $\tilde{I}_\xi$ is uniformly positive definite over $\zeta \in \Xi$, which can be achieved by checking that $\inf_{\zeta \in \Xi} \lambda_{\min}\{P_0(\text{Cov}(Z(\zeta)|V))\} > 0$. We first show that the random vector $(Z, 1\{Z > \zeta\}, Z1\{Z > \zeta\})$ is linearly independent given $V$ pointwisely in $\zeta \in \Xi$. Suppose that given $V$,

$$aZ + b1\{Z > \zeta\} + cZ1\{Z > \zeta\} = 0, \tag{6}$$

a.s., for some constants $a$, $b$ and $c$. Our aim is to show $a = b = c = 0$. When $Z \leq \zeta$, (6) becomes $aZ1\{Z \leq \zeta\} = 0$. Since $\text{Var}(Z|V) > 0$ and $P(Z \leq \zeta|V) > 0$, for every $\zeta \in \Xi$, $\text{Var}(Z|Z \leq \zeta, V) > 0$, and therefore $a = 0$. When $Z > \zeta$, (6) becomes $(b + cZ)1\{Z > \zeta\} = 0$. If $c \neq 0$, $Z = -b/c$, which is contradicted with the fact that $\text{Var}(Z|Z > \zeta, V) > 0$. Thus we conclude that $c = 0$ and $b = 0$ as a consequence. That $P(\text{Cov}(Z(\zeta)|V))$ is uniformly positive definite over $\zeta \in \Xi$ follows since $P(\text{Cov}(Z(\zeta)|V))$ is a continuous function of $\zeta$ and $\Xi$ is compact.

The profile likelihood estimator $\hat{\psi}_n(\zeta)$ can be shown to be consistent for $(\beta_0, \Lambda_0)$ by a similar proof as used for the full maximum likelihood estimator in Huang (1996). The following lemma shows the uniform consistency of $\hat{\psi}_n(\zeta)$ under the null.

LEMMA 3. $\hat{\psi}_n(\zeta) - \psi_0 = o_{P_0}^\Xi(1)$.

To verify the uniform no-bias condition (4), we need the following result about the uniform rate of convergence.

LEMMA 4. Suppose that $d(\eta, \eta_1) : \eta, \eta_1 \in \mathcal{H}_\eta$ is the metric defined on $\mathcal{H}_\eta$, and $C_1$, $C_2$ and $C_3$ are positive constants with,

$$P_0(m_{\beta,\eta,\zeta} - m_{\beta,\eta_0,\zeta}) \leq -C_1 d^2(\eta, \eta_0) + C_2\|\beta - \beta_0\|^2 \tag{7}$$

*and*

$$P_0^\star \sup_{\beta \in B, \eta \in \mathcal{H}_\eta, \|\beta-\beta_0\| < \delta, d(\eta, \eta_0) < \delta, \zeta \in \Xi} |\mathbb{G}_n(m_{\beta, \eta, \zeta} - m_{\beta, \eta_0, \zeta})| \le C_3 \phi_n(\delta), \tag{8}$$

*for functions $\phi_n$ such that $\delta \mapsto \phi_n(\delta)/\delta^\alpha$ is decreasing for some $\alpha < 2$ and sets $B \times \mathcal{H}_\eta \times \Xi$ such that under the null $\Pr(\tilde{\beta}_n \in B, \hat{\eta}_{\tilde{\beta}_n}(\zeta) \in \mathcal{H}_\eta, \zeta \in \Xi) \to 1$. Then $\sup_{\zeta \in \Xi} r_n d(\hat{\eta}_{\tilde{\beta}_n}(\zeta), \eta_0) \le O_{P_0}^\star(1 + r_n\|\tilde{\beta}_n - \beta_0\|)$ for any sequence of positive numbers $r_n$ such that $r_n^2 \phi_n(1/r_n) \le \sqrt{n}$ for every $n$.*

We apply Lemma 4 with $\eta = (\alpha, \Lambda)$, $\mathcal{H}_\eta = \mathbb{R} \times \overline{\mathcal{H}}_\Lambda$, where $\overline{\mathcal{H}}_\Lambda$ is the closed linear span of $\mathcal{H}_\Lambda$, $d(\eta, \eta_1) = \|\alpha - \alpha_1\| + \|\Lambda - \Lambda_1\|_2$ and

$$m_{\beta, \eta, \zeta} = \begin{cases} \log \dfrac{p_{\beta, \eta, \zeta}}{p_{\beta_0, \eta_0}}, & \text{if } \eta = \eta_0, \\[2mm] 2 \log \dfrac{p_{\beta, \eta, \zeta} + p_{\beta_0, \eta_0}}{2 p_{\beta_0, \eta_0}}, & \text{otherwise.} \end{cases}$$

Condition (7) can be established by the Taylor expansion and the uniform boundedness on the derivatives of the loglikelihood. Condition (8) can be verified using Lemma 3.3 of Murphy and van der Vaart (1999), with the choice $\phi_n(\delta) = \delta^{1/2}(1 + M\delta^{-3/2}/\sqrt{n})$, where $M \ge \|m_{\beta, \eta, \zeta}\|_\infty$ is a constant. These conditions imply that $\|\hat{\alpha}_{\tilde{\beta}_n}(\zeta) - \alpha_0\| + \|\hat{\Lambda}_{\tilde{\beta}_n}(\zeta) - \Lambda_0\|_2 = O_P^\Xi(\|\tilde{\beta}_n - \beta_0\| + n^{-1/3})$, for any sequence $\tilde{\beta}_n \to 0$. Now we only need to verify

$$P_0 \dot{\ell}(\beta_0, \beta_0, \hat{\tilde{\eta}}_{\tilde{\beta}_n}(\zeta), \zeta) = o_{P_0}^\Xi(\|\tilde{\beta}_n - \beta_0\| + n^{-1/2}), \tag{9}$$

which is equivalent to (4) under regularity conditions. We further decompose (9) as (17) in Murphy and van der Vaart (2000), which can be easily verified by the Taylor expansion and the uniform boundedness on the first and second derivatives of the loglikelihood.

It is not difficult to see that $\{p_{\xi, \Lambda}(\zeta)\}$ is differentiable in quadratic mean at $(\psi_0, \zeta)$ with respect to the set of directions $\{\xi_0 + th_1, \Lambda_0 + th_2\}$, where $h_1 \in \mathbb{R}^3$, and $h_2$ is a nondecreasing nonnegative right continuous function. Thus all conditions in Section 2 are satisfied.

4.2. *Example 1 revisited: Univariate frailty regression under right censoring.* The odds-rate model we consider in this paper posits that the hazard function has the form (1). We define $g_\zeta(s) \equiv (1 + \zeta s)^{-1/\zeta}$, for $\zeta > 0$, and $g_0(s) \equiv \lim_{\zeta \downarrow 0} g_\zeta(s) = \exp(-s)$. Let $S_Z(\cdot)$ denote the survival function of $T$ given $Z$, and after integrating over $W$, $S_Z(t)$ becomes $g_\zeta(\int_0^t e^{\beta' Z(u)} d\eta(u))$, where the cumulative baseline hazard function $\eta(\cdot)$ is a nonnegative, monotone increasing cadlag (right-continuous with left-hand limits) function. We will argue later that assumptions A–C can be checked under the following conditions. The true null survival function is unique and denoted as $S_0$. The censoring time $C$ is independent of $T$ given $Z$ and uninformative of $\zeta$ and $\beta$. Moreover, for a finite time point $\tau$, $P_0 1\{C \geq \tau\} = P_0 1\{C = \tau\} > 0$ almost surely. $\zeta \in \Xi \equiv [0, K_0]$ for some known $K_0 < \infty$. The null value $\beta_0 = 0$ is an interior point of a known compact set $B_0 \in \mathbb{R}^p$. The parameter space for $\eta$, $\mathcal{H}_\eta$, is a Banach space consisting of continuous and monotone increasing functions on the interval $[0, \tau]$ equipped with the total variation norm $\| \cdot \|_v$. Its closed linear span is denoted as $\overline{\mathcal{H}_\eta}$. The function $\eta(\cdot) \in \mathcal{H}_\eta$ satisfies $\eta(0) = 0$ and $\eta(\tau) < \infty$. The covariate process $Z(\cdot)$ is uniformly bounded in total variation on $[0, \tau]$ and $\mathrm{var}[Z(0+)]$ is positive definite.

The true values of $\pi \equiv (\eta, \zeta)$ are not unique under the null, since the null set $\Pi_0$ contains all pairs of $(\eta, \zeta)$ satisfying, for $t \in [0, \tau]$, $(1 + \zeta\eta(t))^{-1/\zeta} = S_0(t)$, when $\zeta \in (0, K_0]$; and $\exp(-\eta(t)) = S_0(t)$, when $\zeta = 0$. In this example, $\zeta$ appears both under the null and the alternative. Equivalently, for any fixed $\zeta \in (0, K_0]$, $\eta_0(t)(\zeta) = (S_0(t)^{-\zeta} - 1)/\zeta$ and for $\zeta = 0$, $\eta_0(t)(\zeta) = -\log(S_0(t))$, $t \in [0, \tau]$. Hence $\Pi_0 = \{(\zeta, \eta_0(\zeta)) : \zeta \in \Xi\}$. Thus we need a suitable parameter transformation for this example. Let $\overline{\eta} = \phi_\zeta(\eta) \equiv (1 + \eta\zeta)^{1/\zeta} - 1$, for $\zeta > 0$; and $\overline{\eta} = \lim_{\zeta \to 0} \phi_\zeta(\eta) = \exp(\eta) - 1$. It can be easily checked that $\overline{\eta} \in \mathcal{H}_\eta$. The following arguments reveal that the map $\phi_\zeta(\eta) : \mathcal{H}_\eta \mapsto \mathcal{H}_\eta$ is a full-rank reparameterization.

The log likelihood function with the new parameter $\overline{\theta} = (\beta, \overline{\eta}, \zeta)$ is

$$
\begin{aligned}
\ell_n(\overline{\theta}) = \mathbb{P}_n \Big[ & \delta\{\log a_1(v) + (\zeta - 1)\log(\overline{\eta}(v) + 1)\} + \beta' z(v) \\
(10) \qquad & + (1 + \delta\zeta)\log g_\zeta \Big\{ \int_0^v e^{\beta' z(s)} (\overline{\eta}(s) + 1)^{\zeta - 1}\, d\overline{\eta}(s) \Big\} \Big],
\end{aligned}
$$

where $a_1(\cdot)$ is the derivative of $\overline{\eta}(\cdot)$. We will replace $a_1(\cdot)$ with $n\Delta\overline{\eta}(\cdot)$ in the sequel, since this form of the empirical log-likelihood function is asymptotically equal to the true log-likelihood function. When $\beta = 0$, it is clear that $\zeta$ vanishes since (10) $= \mathbb{P}_n\{\delta \log \Delta\overline{\eta}(v) - (\delta + 1)\log(1 + \overline{\eta}(v))\}$, and $\overline{\eta}(0) = 0$. The odds-rate model with new parameterization $\overline{\psi} \equiv (\beta, \overline{\eta})$ is identifiable under the null, since the null survival function $S_0(t|z) = (1 + \overline{\eta})^{-1}$ is a strictly monotone function of $\overline{\eta}$ and is unique.

The Gâteaux derivative of $\phi_\zeta(\eta)$ at $\eta \in \mathcal{H}_\eta$ exists and is obtained by differentiating $\phi_\zeta(\eta)$ along the submodels $t \mapsto \eta + th$. This derivative is $\dot{\phi}_\zeta(\eta)(h) \equiv \partial/\partial t \phi_\zeta(\eta + th)|_{t=0} = (1 + \zeta\eta)^{1/\zeta - 1}h$ for $\zeta > 0$ and $\exp(\eta)h$ for $\zeta = 0$.

The Gâteaux differentiability of $\phi_\zeta(\eta)$ pointwisely in $\zeta$ can be strengthened to uniform Fréchet differentiability by noticing that

$$
\limsup_{t \downarrow 0} \sup_{\zeta \in \Xi} \sup_{\|h\|_v \leq r, h \in \overline{\mathcal{H}_\eta}} \left| \int_0^1 \{\dot{\phi}_\zeta(\eta + sth(\zeta)) - \dot{\phi}_\zeta(\eta)\}\, ds \right| = 0,
$$

for any $r > 0$. Thus $\sup_{\zeta \in \Xi} \sup_{\|h\|_v \leq r, h \in \overline{\mathcal{H}}_\eta} \|\phi_\zeta(\eta + h(\zeta)) - \phi_\zeta(\eta) - \dot{\phi}_\zeta(\eta) \times (h(\zeta))\|_v / \|h(\zeta)\|_v = o(1)$, as $\|h(\zeta)\|_v \to 0$ uniformly over $\zeta \in \Xi$, which we will hereafter refer to as "uniformly Fréchet differentiable." Since $\dot{\phi}_\zeta(\eta)(h)$ is uniformly bounded and Lipschitz in $h$, by checking the definition, we can show that $\dot{\phi}_\zeta$ is one-to-one and continuously invertible uniformly over $\zeta \in \Xi$.

To define a uniformly least-favorable submodel, we calculate scores for $\beta$ and $\overline{\eta}$. Let $\mathcal{H}$ denote the space of elements $h = (h_1, h_2)$ such that $h_1 \in \mathbb{R}^p$ and $h_2 \in \mathcal{H}_\eta$. Consider the one-dimensional submodel defined by the map $t \mapsto \overline{\psi}_t \equiv \overline{\psi} + t(h_1, \int_0^{(\cdot)} h_2(u) \, d\overline{\eta}(u))$, $h \in \mathcal{H}$. The derivative of $\log \ell_n(\overline{\psi}_t, \zeta)$ with respect to $t$ evaluated at $t = 0$ yields score operators $\dot{\ell}_n(\overline{\psi}, \zeta)(h) \equiv (\dot{\ell}_{n\beta}(h_1), \dot{\ell}_{n\overline{\eta}}(h_2))$, where

$$\dot{\ell}_{n\beta}(\overline{\psi}, \zeta)(h_1)$$
$$= \mathbb{P}_n \dot{\ell}_\beta(h_1)$$
$$= \mathbb{P}_n \bigg\{ \delta h_1' Z(X)$$
$$\qquad - (1 + \delta\zeta) \frac{\int_0^\tau h_1' Z(u) Y(u) e^{\beta' Z(u)} (\overline{\eta}(u) + 1)^{\zeta - 1} \, d\overline{\eta}(u)}{1 + \zeta \int_0^\tau h_1' Z(u) Y(u) e^{\beta' Z(u)} (\overline{\eta}(u) + 1)^{\zeta - 1} \, d\overline{\eta}(u)} \bigg\},$$

and

$$\dot{\ell}_{n\overline{\eta}}(\overline{\psi}, \zeta)(h_2)$$
$$= \mathbb{P}_n \dot{\ell}_{\overline{\eta}}(h_2)$$
$$= \mathbb{P}_n \bigg\{ \int_0^\tau \bigg( h_2(u) + \frac{(\zeta - 1) \int_0^u h_2(s) d\overline{\eta}(s)}{\overline{\eta}(u) + 1} \bigg) \, dN(u)$$
$$\qquad - (1 + \delta\zeta)$$
$$\qquad\qquad \times \int_0^\tau Y(u) e^{\beta' Z(u)} (\overline{\eta}(u) + 1)^{\zeta - 2}$$
$$\qquad\qquad\qquad \times \bigg[ (\zeta - 1) \int_0^u h_2(s) \, d\overline{\eta}(s) + h_2(u)(1 + \overline{\eta}(u)) \bigg] d\overline{\eta}(u)$$
$$\qquad\qquad\qquad \times \bigg( 1 + \zeta \int_0^\tau Y(u) e^{\beta' Z(u)} (\overline{\eta}(u) + 1)^{\zeta - 1} \, d\overline{\eta}(u) \bigg)^{-1} \bigg\},$$

with $Y(u) \equiv 1\{V \geq u\}$.

To obtain the information operator, we consider the two-dimensional submodel defined by the map $(s, t) \mapsto \overline{\psi}_{st} \equiv \overline{\psi} + s(h_1, \int_0^{(\cdot)} h_2(u) \, d\overline{\eta}(u)) + t(\tilde{h}_1, \int_0^{(\cdot)} \tilde{h}_2(u) \, d\overline{\eta}(u))$, where $h, \tilde{h} \in \mathcal{H}$. Define $\mathcal{H}_\infty = \{h \in \mathcal{H} : \|h\|_{\mathcal{H}} < \infty\}$. The information operator $\overline{\sigma}_{\overline{\eta}}(h) : \mathcal{H}_\infty \mapsto \mathcal{H}_\infty$ is given by $-P_0 \partial / \partial s \partial t \ell_1(\overline{\psi}_{st})|_{s,t=0} = \overline{\psi}(\overline{\sigma}_{\overline{\eta}}(h))$. We will show $\overline{\sigma}_{\overline{\eta}}$ is one-to-one, continuously invertible and onto uniformly over $\zeta \in \Xi$, via Part (1) of Lemma 7 in the Section 7 for which it suffices to show that the information operator for the original parameterization $\sigma_\theta$ is one-to-one, continuously invertible and onto uniformly over $\zeta \in \Xi$.

With the same derivation of $\overline{\sigma}_{\overline{\theta}}$, $\sigma_\theta : \mathcal{H}_\infty \mapsto \mathcal{H}_\infty$ takes the form

$$\sigma_\theta(h) = \begin{pmatrix} \sigma_\theta^{11} & \sigma_\theta^{12} \\ \sigma_\theta^{21} & \sigma_\theta^{22} \end{pmatrix} \begin{pmatrix} h_1 \\ h_2 \end{pmatrix},$$

where

$$\sigma_\theta^{11}(h_1) = -P_0 S(\theta) \int_0^\tau h_1' Z(u) Y(u) e^{\beta' Z(u)} \, d\eta_0(u),$$

$$\sigma_\theta^{12}(h_2) = -P_0 S(\theta) \int_0^\tau h_2(u) Z(u) Y(u) e^{\beta' Z(u)} \, d\eta_0(u),$$

$$\sigma_\theta^{21}(h_1) = -P_0 S(\theta)(1 + \zeta\eta(T \wedge \tau)(\zeta)) h_1' Z(u) Y(u) e^{\beta' Z(u)}$$
$$- \zeta P_0 S(\theta) Y(u) \int_0^\tau h_1' Z(u) Y(u) e^{\beta' Z(u)} \, d\eta_0(u),$$

$$\sigma_\theta^{22}(h_2) = -P_0 S(\theta)(1 + \zeta\eta(T \wedge \tau)(\zeta)) h_2(u) Y(u) e^{\beta' Z(u)}$$
$$- \zeta P_0 S(\theta) Y(u) \int_0^\tau h_2(u) Y(u) e^{\beta' Z(u)} \, d\eta_0(u),$$

with $S(\theta) = -(1 + \delta\zeta)/(1 + \zeta\overline{\eta}(\tau))^2$.

All of the operators $\sigma_\theta^{ij}$, $1 \leq i, j \leq 2$ are uniformly compact and bounded over $\zeta \in \Xi$. With a similar argument as in Kosorok, Lee and Fine (2004), the linear operator $\sigma_\theta : \mathcal{H}_\infty \mapsto \mathcal{H}_\infty$ is one-to-one, continuously invertible and onto uniformly over $\zeta \in \Xi$ by verifying the conditions of Lemma 8 in the Section 7. Thus a uniformly least-favorable submodel for estimating $\beta$ in the presence of $\overline{\eta}$ and $\zeta$ is $\overline{\eta}_t(\beta, \overline{\eta}, \zeta) = (1 + (\beta - t)'\nu_{\overline{\eta}}) \, d\overline{\eta}$, where $\nu_{\overline{\eta}} : \mathbb{R} \mapsto \mathbb{R}^p$ is the uniformly least-favorable direction at $(\beta_0, \overline{\eta}, \zeta)$ defined by $h'\nu_{\overline{\eta}} = (\overline{\sigma}_{\overline{\theta}}^{22})^{-1} \overline{\sigma}_{\overline{\theta}}^{21} h$, $h \in \mathbb{R}^p$. This leads to $\ell(t, \beta, \overline{\eta}, \zeta) = \ell_1(\beta, \overline{\eta}_t(\theta), \zeta)$. Because $\overline{\eta}_{\overline{\eta}}(\beta, \overline{\eta}, \zeta) = \overline{\eta}$, B1 is satisfied. Since $\partial/\partial t|_{t=\beta_0} \ell(t, \beta_0, \overline{\eta}_0, \zeta) = \dot{\ell}_\beta(\beta_0, \overline{\psi}_0, \zeta) = \tilde{\ell}_\beta(\overline{\psi}_0, \zeta)$, where $\tilde{\ell}_\beta(x) = \dot{\ell}_\beta - \dot{\ell}_{\overline{\eta}}\nu_\theta$ is the efficient score for $\beta$, B2 is satisfied due to the continuity of the involved functions with respect to $\overline{\psi}$ and the fact that $\Xi$ is compact. The efficient information for $\beta$ is $\tilde{\mathcal{I}}_\beta = P_0 \tilde{\ell}_\beta \tilde{\ell}_\beta'$. That $\{\dot{\ell}(t, \overline{\psi}, \zeta) : t \in U, \overline{\psi} \in V, \zeta \in \Xi\}$ is $P_0$-Donsker and $\{\ddot{\ell}(t, \overline{\psi}, \zeta) : t \in U, \overline{\psi} \in V, \zeta \in \Xi\}$ is $P_0$-Glivenko–Cantelli for some neighborhoods $U$ and $V$ follows from standard empirical process arguments.

It follows from Corollary 8.1.3 in Golub and Van Loan (1983) that the set of eigenvalues is a continuous function of the elements of $\tilde{\mathcal{I}}_\beta(\overline{\theta})$, which are continuous functions of $\zeta$. The set of eigenvalues is therefore a continuous function of $\zeta$. Thus $\inf_\zeta \{\lambda_{\min}\{\tilde{I}_\beta(\theta_0(\zeta))\}\} > 0$ by the compactness of $\Xi$, and assumption A1 is satisfied.

The consistency of the restricted MLE $\hat{\tilde{\psi}}_0$ and the uniform consistency of the unrestricted MLE $\hat{\tilde{\psi}}_n(\zeta)$ can be established via the self-consistency equation approach, with arguments similar to the proof of Theorem 3 in Kosorok, Lee and Fine (2004). We omit the details. To verify the uniform no-bias condition (4), it suffices to show that

$$\sup_{\zeta \in \Xi} \| \hat{\overline{\eta}}_{\tilde{\beta}_n}(\zeta) - \overline{\eta}_0 \|_\infty = O^\star_{P_0}(\| \tilde{\beta}_n - \beta_0 \| + n^{-1/2})$$

for any sequence $\tilde{\beta}_n \to \beta_0$,

where "$\star$" denotes outer probability. By verifying conditions in Lemma 9 in the Section 7, we have

$$\sup_{\zeta \in \Xi}(\mathbb{P}_n - P_0)\{\dot{\ell}_{\overline{\psi}}(\tilde{\beta}_n, \hat{\overline{\eta}}_{\tilde{\beta}_n}(\zeta), \zeta) - \dot{\ell}_{\overline{\psi}}(\beta_0, \overline{\eta}_0, \zeta)\} = o^\star_{P_0}(n^{-1/2}).$$

Together with the fact that $\mathbb{P}_n \dot{\ell}_{\overline{\psi}}(\tilde{\beta}_n, \hat{\overline{\eta}}_{\tilde{\beta}_n}(\zeta), \zeta) = P_0 \dot{\ell}_{\overline{\psi}}(\beta_0, \overline{\eta}_0, \zeta) = 0$, we obtain

$$\begin{aligned}
P_0 &\{\dot{\ell}_{\overline{\psi}}(\tilde{\beta}_n, \hat{\overline{\eta}}_{\tilde{\beta}_n}(\zeta), \zeta) - \dot{\ell}_{\overline{\psi}}(\beta_0, \overline{\eta}_0, \zeta)\} \\
&= P_0 \dot{\ell}_{\overline{\psi}}(\tilde{\beta}_n, \hat{\overline{\eta}}_{\tilde{\beta}_n}(\zeta), \zeta) - \mathbb{P}_n \dot{\ell}_{\overline{\psi}}(\tilde{\beta}_n, \hat{\overline{\eta}}_{\tilde{\beta}_n}(\zeta), \zeta) \\
&= -(\mathbb{P}_n - P_0)\dot{\ell}_{\overline{\psi}}(\beta_0, \overline{\eta}_0, \zeta) + o^\star_{P_0}(n^{-1/2}),
\end{aligned}$$

uniformly over $\zeta \in \Xi$.

Let $\dot{l}_\psi(h) \equiv (\dot{l}_\beta(h_1), \dot{l}_\eta(h_2))$ denote the score operator of $\psi$ with the original parameterization. It was shown in Kosorok, Lee and Fine (2004) that the operator $\psi \mapsto \dot{l}_\psi$ is Fréchet differentiable with derivative $\psi(\sigma_\theta(h))$, and it can be strengthened to uniform Fréchet differentiability due to the smoothness of the involved functions. Since $\phi_\zeta$ is uniformly Fréchet differentiable, by Part (2) of Lemma 7, the chain rule for uniform Fréchet differentiability, $\dot{\ell}_{\overline{\psi}} \equiv (\dot{\ell}_\beta, \dot{\ell}_{\overline{\eta}})$ is uniformly Fréchet differentiable with derivative $\sigma_{\phi_\zeta^{-1}(\overline{\theta})} \circ \dot{\phi}_\zeta^{-1}(\overline{\theta})$.

By the uniform Fréchet differentiability of $\dot{\ell}_{\overline{\psi}}$,

$$\begin{aligned}
\overline{\sigma}_{\overline{\theta}}(\tilde{\beta}_n, \hat{\overline{\eta}}_{\tilde{\beta}_n}(\zeta) - \overline{\eta}_0) &= P_0\{\dot{\ell}_{\overline{\psi}}(\tilde{\beta}_n, \hat{\overline{\eta}}_{\tilde{\beta}_n}(\zeta), \zeta) - \dot{\ell}_{\overline{\psi}}(\beta_0, \overline{\eta}_0, \zeta)\} \\
&\quad + o^{\overline{\Xi}}_{P_0}(\| \tilde{\beta}_n - \beta_0 \| + \| \hat{\overline{\eta}}_{\tilde{\beta}_n}(\zeta) - \overline{\eta}_0 \|_\infty).
\end{aligned}$$

Since $\overline{\sigma}_{\overline{\theta}}$ is linear, the first term on the right-hand side is of the order $O_{P_0}(n^{-1/2})$. It follows that $\sup_{\zeta \in \Xi} \| \hat{\overline{\eta}}_{\tilde{\beta}_n}(\zeta) - \overline{\eta}_0 \|_\infty = O^\star_{P_0}(\| \tilde{\beta}_n - \beta_0 \| + n^{-1/2})$, since $\overline{\sigma}_{\overline{\theta}}$ is uniformly continuously invertible over $\zeta \in \Xi$.

**5. Simulation results.** This section presents simulation results regarding the finite-sample properties of the proposed optimal test statistics for Example 2, the change-point Cox model with current status data. The simulation study was designed with several objectives. First, we demonstrate how to compute the asymptotic critical values with the proposed weighted bootstrap procedure. Second, we analyze the empirical type I error of the proposed tests and compare with the nominal size of the tests. Third, we compare the power of the optimal tests with that of other tests such as the sup score statistics (equivalent to the likelihood ratio statistic) and some naive (pointwise) tests under several different alternatives. Fourth, we evaluate the sensitivity of the power of the optimal test to the choice of $c$ under several different alternatives.

A single time-independent covariate $Z$ with a uniform $[0,1]$ distribution was used. The threshold covariate $Y = Z$. The parameter $\alpha$ was set at $\alpha_0 = 0$, with the cumulative baseline hazards $A_0(t) = 3t^2$. The censoring time was uniformly distributed on the interval $[0,5]$. This resulted in a censoring rate of about 25% under the null hypothesis. Under the alternative, we set $\beta_{10} = -0.5$, $\beta_{20} \in \{-0.3, -0.5, -0.8\}$. The range of $\beta_{20}$ values reflects the distance from the null. We consider the following alternative distributions of $\zeta$:

1. The weight $J(\zeta)$ degenerates to one point at 0.5, that is $\zeta = 0.5$.
2. A uniform weight $J(\zeta)$ with support on $[0.05, 0.95]$.

For all the scenarios, we compute the optimal tests with a uniform weight on $[0.05, 0.95]$. The sample size for each simulated data set was 300. For each simulated data set, 250 bootstraps were generated with standard exponential weights truncated at 5, to compute the critical values for $R_n(\zeta)$, the naive score statistic at several $\zeta$ values, $\sup_\zeta R_n(\zeta)$, the sup score statistic and $ER_n$, the weighted exponential score statistic. We take $c = 0, 0.5, 1, 3$ and $\infty$, respectively. Each scenario was replicated 1000 times. To compute the restricted MLE under the null, we use the iterative convex minorant algorithm. Empirical type I error and power results for selected subsets of the test statistics described above are provided in Table 1.

We now make several general comments on the simulation results. The empirical type I error for all the tests is quite close to the nominal level. When the alternative distribution of $\zeta$ is correctly specified, the optimal test is notably more powerful than the sup score statistic and naive tests. When the true alternative distribution of $\zeta$ degenerates to one point, although the weighted exponential tests are no longer optimal, the empirical powers are still superior to the naive tests with misspecified $\zeta$. We also observe that the empirical power of the sup score statistics is comparable to that of the naive test at the true $\zeta$, which may be due to the fast convergence rate of the change-point estimator. For all the alternatives considered, the empirical

TABLE 1
*The empirical type* I *error and power of the proposed tests, sample size* $n = 300$, *1000 simulations, with bootstrap size 250. The worst case Monte Carlo error for table entries is 0.016. The Monte Carlo error is 0.007 and 0.009 for empirical type* I *error with nominal size 0.05 and 0.1, respectively. The empirical power results are based on size 0.05 tests*

**Empirical type I error**

| Nominal size | Weighted exponential tests, $c =$ | | | | | Sup score | Naive tests $R_n(\zeta)$, $\zeta =$ | | |
|---|---|---|---|---|---|---|---|---|---|
| | **0** | **0.5** | **1** | **3** | **∞** | | **0.3** | **0.9** | **0.5** |
| $J(\zeta) \sim \text{Uniform}[0.05, 0.95]$ | | | | | | | | | |
| 0.05 | 0.056 | 0.057 | 0.045 | 0.063 | 0.046 | 0.058 | 0.043 | 0.044 | 0.039 |
| 0.10 | 0.098 | 0.103 | 0.109 | 0.095 | 0.099 | 0.085 | 0.112 | 0.103 | 0.100 |

**Empirical power**

| True alternative | Weighted exponential tests, $c =$ | | | | | Sup score | Naive tests $R_n(\zeta)$, $\zeta =$ | | |
|---|---|---|---|---|---|---|---|---|---|
| | **0** | **0.5** | **1** | **3** | **∞** | | **0.3** | **0.9** | **0.5** |
| $J(\zeta) \sim \text{Uniform}[0.05, 0.95]$ | | | | | | | | | |
| $\zeta = 0.5$ | | | | | | | | | |
| $\eta = -0.3$ | 0.646 | 0.647 | 0.653 | 0.653 | 0.656 | 0.688 | 0.243 | 0.044 | 0.692 |
| $\eta = -0.5$ | 0.835 | 0.833 | 0.839 | 0.845 | 0.847 | 0.865 | 0.616 | 0.076 | 0.840 |
| $\eta = -0.8$ | 0.922 | 0.925 | 0.928 | 0.928 | 0.928 | 0.968 | 0.957 | 0.174 | 0.942 |
| $J(\zeta) \sim \text{Uniform}[0.05, 0.95]$ | | | | | | | | | |
| $\eta = -0.3$ | 0.320 | 0.320 | 0.320 | 0.320 | 0.312 | 0.211 | 0.133 | 0.055 | 0.142 |
| $\eta = -0.5$ | 0.485 | 0.488 | 0.492 | 0.494 | 0.500 | 0.405 | 0.258 | 0.083 | 0.272 |
| $\eta = -0.8$ | 0.748 | 0.757 | 0.763 | 0.768 | 0.769 | 0.605 | 0.494 | 0.183 | 0.413 |

power of the weighted exponential tests seems to increase as $c$ increases. However, the trend is rather weak. In many cases, the difference in power is less than 0.01. This suggests that the direction of the test (specifically, least favorable curve in this paper), rather than the scale of the curve, is most critical for the power of the weighted exponential test.

**6. Discussion.** In this paper, we consider tests of hypotheses when the parameters are not identifiable under the null in semiparametric models. Our optimality results apply to a large class of semiparametric testing problems under loss of identifiability, where nuisance parameters may not be root-$n$ estimable either under the null or alternative. We note that our current regularity conditions are not directly applicable for testing under loss of inevitability when the parameter of interest is not root-$n$ estimable. One example is testing for homogeneity in mixture models, where the usual first order Taylor approximation may not be possible [Chen, Chen and Kalbfleisch (2004); Chernoff and Lander (1995); Dacunha-Castelle and Gassiat (1999); Lindsay (1995); Liu and Shao (2003)]. A higher order expansion is required.

Although not directly covered by our framework, the homogeneity tests may possess a uniform quadratic expansion [Zhu and Zhang (2006)], thus permitting a generalization of our results to general quadratic expansions. In the following, we conclude the paper with a brief discussion of this generalization.

To be concrete, let us consider a two-component mixture with density $g(\rho, \mu_1, \mu_2, \eta) = \rho f(\mu_1, \eta) + (1 - \rho) f(\mu_2, \eta)$, where $f(\mu, \eta)$ is a parametric p.d.f. with parameters $\mu \in \mathbb{R}^p$, $\eta \in \mathbb{R}^q$, such as a location-scale family. Let $\beta = \mu_2 - \mu_1$, $\theta = (\rho, \beta', \mu_1', \eta')'$, and the hypothesis of interest is $\beta = 0$; that is, there is only a single component in the mixture. For convenience, we assume the mixing proportion $\rho \in (0, 1]$ and $\mu_1 = \mu_2 = \mu_0$ under the null.

In this example, $\rho$ is not identifiable and $\mu_1$ and $\mu_2$ are mutually indistinguishable under the null. Simple algebra shows that the information matrix for $\psi \equiv (\beta, \mu_1)$ is singular under the null, for arbitrary values of $\rho$, which corresponds to the fact that $\mu_1$ and $\mu_2$ are not root-$n$ estimable [Chen and Chen (2003); Zhu and Zhang (2004)]. We consider the following reparameterization: $\overline{\mu} = (1 - \rho)(\mu_1 - \mu_0) + \rho(\mu_2 - \mu_0)$ and $\overline{v} = (1 - \rho)(\mu_1 - \mu_0)^2 + \rho(\mu_2 - \mu_0)^2$, which can be considered as "mixed mean" and "mixed variance". Let $\beta \equiv (\overline{\mu}, \overline{v})$ and $\overline{\psi} \equiv (\overline{\mu}, \overline{v}, \eta)$. We can establish the identifiability of $\overline{\psi}$ and the consistency and the root-$n$ rate of the MLE of $\overline{\psi}$ under the null. Furthermore, under a set of assumptions on the parameter space [e.g., the cone condition in Andrews (1999, 2001)] and the stochastic differentiability and equicontinuity of the involved functions, we can establish the following quadratic expansion of the loglikelihood with respect to $\overline{\psi}$:

$$L_n(\overline{\psi}, \zeta) = L_n(\overline{\psi}_0, \zeta) + (\overline{\psi} - \overline{\psi}_0)' S_{n\zeta}(\overline{\psi}_0)$$
$$+ \tfrac{1}{2}(\overline{\psi} - \overline{\psi}_0)' B_\zeta(\overline{\psi}_0)(\overline{\psi} - \overline{\psi}_0) + r_n(\zeta),$$

where $r_n(\zeta) = o_p^{\overline{\Xi}}(1)$, and $S_{n\zeta}(\cdot)$ and $B_\zeta(\cdot)$ are different from but similar in structure to the score and information processes for $\overline{\psi}$ indexed by $\zeta$.

When the nuisance parameter $\eta$ is not present, a similar weight as in the current paper for $\overline{\psi}$ can be chosen as $Q_\zeta(\cdot) = B_\zeta(\cdot)$. The corresponding weighted exponential tests are still optimal in the Neyman–Pearson sense. If $\eta$ is present, a uniformly least favorable curve for this quadratic expansion with respect to $\beta$ would need to be characterized. This is beyond the scope of the current paper but is an interesting topic for future research.

## 7. Proofs.

PROOF OF LEMMA 1.   Since $\dot{\phi}_\zeta$ is linear, continuously invertible and one-to-one, the tangent set for $\eta$ and $\overline{\eta}$ are identical. By the chain rule, $\dot{\ell}_{\overline{\eta}}(\gamma) = \dot{\ell}_\eta \dot{\phi}_\zeta^{-1}(\gamma)$ for any $\gamma$ in the tangent set of $\overline{\eta}$. The efficient score for $\beta$ with the parameter $(\beta, \eta, \zeta)$ is: $\tilde{\ell}_\beta(\beta, \eta, \zeta) = (I - \dot{\ell}_\eta(\dot{\ell}_\eta^\star \dot{\ell}_\eta)^{-1} \dot{\ell}_\eta^\star) \dot{\ell}_\beta(\psi, \zeta)$ and with the parameter $(\beta, \overline{\eta}, \zeta)$ is: $(I - \dot{\ell}_{\overline{\eta}}(\dot{\ell}_{\overline{\eta}}^{\star}\dot{\ell}_{\overline{\eta}})^{-1}\dot{\ell}_{\overline{\eta}}^{\star})\dot{\ell}_{\beta}(\psi, \zeta)$. The efficient score function is invariant under such reparameterizations since

$$
\begin{aligned}
I - \dot{\ell}_{\overline{\eta}}(\dot{\ell}_{\overline{\eta}}^{\star}\ell_{\overline{\eta}})^{-1}\dot{\ell}_{\overline{\eta}}^{\star}(\overline{\psi}, \zeta) &= I - \dot{l}_{\eta}\dot{\phi}_{\zeta}^{-1}(\dot{\phi}_{\zeta}^{-1\star}\dot{l}_{\eta}^{\star}\dot{l}_{\eta}\dot{\phi}_{\zeta}^{-1})^{-1}\dot{\phi}_{\zeta}^{\star}\dot{l}_{\eta}^{\star}(\psi, \zeta) \\
&= I - \dot{l}_{\eta}\dot{\phi}_{\zeta}^{-1}\dot{\phi}_{\zeta}(\dot{l}_{\eta}^{\star}\dot{l}_{\eta})^{-1}\dot{\phi}_{\zeta}^{\star}(\dot{\phi}_{\zeta}^{\star})^{-1}\dot{l}_{\eta}^{\star}(\psi, \zeta) \\
&= I - \dot{l}_{\eta}(\dot{l}_{\eta}^{\star}\dot{l}_{\eta})^{-1}\dot{l}_{\eta}^{\star}(\psi, \zeta),
\end{aligned}
$$

and $\dot{\ell}_{\beta}(\overline{\psi}, \zeta) = \dot{l}_{\beta}(\psi, \zeta)$. That the efficient information matrix is invariant under reparameterizations thus follows from its definition. $\square$

PROOF OF LEMMA 2. It suffices to show that under the full rank reparameterization, for any random sequence $\tilde{\beta}_n \to_{P_0} \beta_0$,

$$
\begin{aligned}
(11) \quad \log pl_n(\tilde{\beta}_n, \zeta) = {} &\log pl_n(\beta_0, \zeta) + n(\tilde{\beta}_n - \beta_0)'\mathbb{P}_n\tilde{\ell}_{\beta}(\overline{\psi}_0, \zeta) \\
&- \tfrac{1}{2}n(\tilde{\beta}_n - \beta_0)'\tilde{\mathcal{I}}_{\beta}(\overline{\psi}_0, \zeta)(\hat{\beta}_n - \beta_0) \\
&+ o_{P_0}^{\overline{\overline{=}}}(\sqrt{n}\|\tilde{\beta}_n - \beta_0\| + 1)^2.
\end{aligned}
$$

By assumptions B2, B4 and the dominated convergence theorem, for every $(\tilde{t}, \tilde{\beta}, \tilde{\overline{\eta}}) - (\beta_0, \beta_0, \overline{\eta}_0) \to 0$, we have $P_0(\dot{\ell}(\tilde{t}, \tilde{\beta}, \tilde{\overline{\eta}}, \zeta) - \dot{\ell}_{\beta}(\overline{\psi}_0, \zeta))^2 = o^{\overline{=}}(1)$. Similarly, we have $P_0\ddot{\ell}(\tilde{t}, \tilde{\beta}, \tilde{\overline{\eta}}, \zeta) - P_0\ddot{\ell}(\beta_0, \beta_0, \overline{\eta}_0, \zeta) = o^{\overline{=}}(1)$. The derivative of the function $t \mapsto \log \ell(t, \overline{\psi}_0, \zeta)$ satisfies $P_0\ddot{\ell}(\beta_0, \overline{\psi}_0, \zeta) = -\tilde{\mathcal{I}}_{\beta}(\overline{\psi}_0, \zeta)$. These facts, together with the empirical process conditions, imply that for every random sequence $(\tilde{t}, \tilde{\beta}, \tilde{\overline{\eta}}) \to (\beta_0, \beta_0, \overline{\eta}_0)$, $\mathbb{G}_n\dot{\ell}(\tilde{t}, \tilde{\beta}, \tilde{\overline{\eta}}) - \mathbb{G}_n\tilde{\ell}_{\beta}(\overline{\psi}_0, \zeta) = o_{P_0}^{\overline{\overline{=}}}(1)$ and $\mathbb{P}_n\ddot{\ell}(\tilde{t}, \tilde{\beta}, \tilde{\overline{\eta}}, \zeta) + \tilde{\mathcal{I}}_{\beta}(\overline{\psi}_0, \zeta) = o_{P_0}^{\overline{\overline{=}}}(1)$. The subsequent steps of the proof are similar to those used in the proof of Theorem 1 in Murphy and van der Vaart (2000), and we omit the details. $\square$

PROOF OF THEOREM 1. The proof takes several steps. We first show the asymptotic equivalence of these statistics, which is summarized in Lemma 5 below. With a small abuse of notation, let $PLR_n \equiv \int pl_n(\beta + h/\sqrt{n}, \zeta)\,dQ_{\zeta}(h)\,dJ(\zeta)/$
$pl_n(\beta_0, \zeta)$. This is the profile likelihood ratio of the alternative over the null and it can be approximated by

$$
\begin{aligned}
\overline{PLR}_n \equiv {} &\int \exp\{\tfrac{1}{2}\overline{\beta}_n(\theta_0(\zeta))'\tilde{I}_{\beta}(\theta_0(\zeta))\overline{\beta}_n(\theta_0(\zeta))\} \\
&\times \int \exp\{-\tfrac{1}{2}(\overline{\beta}_n(\theta_0(\zeta)) - h)' \\
&\qquad \times \tilde{I}_{\beta}(\theta_0(\zeta))(\overline{\beta}_n(\theta_0(\zeta)) - h)\}\,dQ_{\zeta}(h)\,dJ(\zeta),
\end{aligned}
$$

with the linear statistic $\overline{\beta}_n(\theta_0(\zeta)) \equiv \sqrt{n}\tilde{I}_\beta^{-1}(\theta_0(\zeta))\mathbb{P}_n\tilde{l}_\beta(\theta_0(\zeta))$. An approximate exponential Wald statistic $\overline{EW}_n$ is defined as

$$\overline{EW}_n = (1+c)^{-p/2}\int \exp\left(\frac{1}{2}\frac{c}{1+c}\overline{W}_n(\zeta)\right)dJ(\zeta),$$

where $\overline{W}_n(\zeta) = \overline{\beta}_n(\theta_0(\zeta))'\tilde{I}_\beta(\theta_0(\zeta))\overline{\beta}_n(\theta_0(\zeta))$.

Now we show the asymptotic distribution of these tests under the null hypothesis. Assume without loss of generality that $\hat{\beta}_n$ and $\hat{\overline{\psi}}_n$ take their values in $U$ and $V$ as defined in assumption B4, respectively. Following Lemma 3.2 in Murphy and van der Vaart (1997), we have $\mathbb{G}_n(\hat{\ell}(\hat{\beta}_n, \hat{\overline{\psi}}_n, \zeta) - \tilde{\ell}_\beta(\overline{\psi}_0, \zeta)) \to_{P_0} 0$. Thus $\tilde{\ell}_\beta(\psi_0, \zeta) = \tilde{l}_\beta(\theta_0(\zeta))$ is $P_0$-Donsker as a class indexed by $\zeta \in \Xi$ and $\overline{EW}_n \to_d e\chi(c)$ by the continuous mapping theorem. Lemma 5 below then gives the desired results of Theorem 1. □

LEMMA 5. *Under the null hypothesis and assumptions* A–C, (1) $PLR_n - \overline{PLR}_n \to_{P_0} 0$, (2) $\overline{PLR}_n = \overline{EW}_n$, (3) $\overline{EW}_n - EW_n \to_{P_0} 0$, (4) $\overline{EW}_n - ER_n \to_{P_0} 0$ *and* (5) $ER_n - ELR_n \to_{P_0} 0$.

PROOF. For notational simplicity, let $\overline{\beta}_n = \overline{\beta}_n(\theta_0(\zeta))$ and $\tilde{I}_0 = \tilde{I}_\beta(\theta_0(\zeta))$. We first show (1). For $0 < M < \infty$, define

$$PLR_n(M) = \int_{\zeta \in \Xi}\int_{\|h\| \le M} pl_n(\beta_0 + h/\sqrt{n}, \zeta)\,dQ_\zeta(h)\,dJ(\zeta)/pl_n(\beta_0, \zeta),$$

and

$$\overline{PLR}_n(M) = \int_{\zeta \in \Xi}\exp(\tfrac{1}{2}\overline{\beta}_n'\tilde{I}_\beta\overline{\beta}_n)$$
$$\times \int_{\|h\| \le M}\exp(-\tfrac{1}{2}(\overline{\beta}_n - h)'\tilde{I}_0(\overline{\beta}_n - h))\,dQ_\zeta(h)\,dJ(\zeta).$$

Note that for any $M > 0$,

$$|PLR_n - \overline{PLR}_n| \le |PLR_n - PLR_n(M)| + |PLR_n(M) - \overline{PLR}_n(M)| + |\overline{PLR}_n - \overline{PLR}_n(M)|.$$

Hence it suffices to show that (i) $|PLR_n - PLR_n(M)| \to_{P_0} 0$, (ii) $|\overline{PLR}_n - \overline{PLR}_n(M)| \to_{P_0} 0$ and (iii) $|PLR_n(M) - \overline{PLR}_n(M)| \to_{P_0} 0$, as $n \to \infty$ and $\forall M : 0 < M < \infty$. To show (i), for any $\varepsilon > 0$,

$$\Pr(|PLR_n - PLR_n(M)| > \varepsilon)$$
$$\le \varepsilon^{-1}P_0|PLR_n - PLR_n(M)|$$
$$(12) \qquad = \varepsilon^{-1}P\int_{\zeta \in \Xi}\int_{\|h\| > M}\frac{pl_n(\beta_0 + h/\sqrt{n}, \zeta)}{pl_n(\beta_0, \zeta)}\,dQ_\zeta(h)\,dJ(\zeta)$$

$$(13) \qquad \leq \varepsilon^{-1} P \int_{\zeta \in \Xi} \int_{\|h\| > M} \frac{pl_n(\hat{\beta}_n(\zeta), \zeta)}{pl_n(\beta_0, \zeta)} \, dQ_\zeta(h) \, dJ(\zeta)$$

$$(14) \qquad \to \varepsilon^{-1} P \int_{\zeta \in \Xi} \int_{\|h\| > M} (1 + o_p(1)) \, dQ_\zeta(h) \, dJ(\zeta)$$

$$(15) \qquad = \varepsilon^{-1} \int_{\zeta \in \Xi} \int_{\|h\| > M} dQ_\zeta(h) \, dJ(\zeta) + o(1),$$

where (12) uses assumption C and (13) holds by definition of the profile likelihood. (14) holds by assumption B3 and Lemma 2. (15) holds by Fubini's theorem. The right-hand side of (15) can be made arbitrarily small for all $n$ by taking $M$ large enough, since $Q_\zeta$ is a uniformly tight measure.

For (ii), we have

$$
\begin{aligned}
(16) \quad & |\overline{PLR}_n - \overline{PLR}_n(M)| \\
&= \int_{\zeta \in \Xi} \exp(\tfrac{1}{2} \mathbb{P}_n \tilde{l}_\beta(\theta_0(\zeta))' \tilde{I}_0^{-1} \mathbb{P}_n \tilde{l}_\beta(\theta_0(\zeta))) \\
&\quad \times \int_{\|h\| > M} \exp(-\tfrac{1}{2}(\overline{\beta}_n - h)' \tilde{I}_0(\overline{\beta}_n - h)) \, dQ_\zeta(h) \, dJ(\zeta) \\
&\leq \exp\left(\tfrac{1}{2} \sup_{\zeta \in \Xi} \|\mathbb{P}_n \tilde{l}_\beta(\theta_0(\zeta))\|^2 \sup_{\zeta \in \Xi} \|\tilde{I}_0^{-1}\|\right) \int_{\zeta \in \Xi} \int_{\|h\| > M} dQ_\zeta(h) \, dJ(\zeta).
\end{aligned}
$$

In the inequality, $\|\mathbb{P}_n \tilde{l}_\beta(\theta_0(\zeta))\|^2 = O_{P_0}^{\overline{\Xi}}(1)$ follows from assumption B4. The fact that $\|\tilde{I}_0^{-1}\| = O_{P_0}^{\overline{\Xi}}(1)$ follows from assumption A. The last term $\int_{\zeta \in \Xi} \int_{\|h\| > M} dQ_\zeta(h) \, dJ(\zeta) \to 0$, as $M \to \infty$. Hence (16) $= o_p(1)$, as $M \to \infty$.

Now we show (iii). For contiguous sequences $\beta_0 + h/\sqrt{n} \to_{P_0} \beta_0$ and $\|h\| \leq M$, Lemma 2 yields the following expansion of the profile likelihood under the null:

$$
\begin{aligned}
\log pl_n(\beta_0 + h/\sqrt{n}, \zeta) &= \log pl_n(\beta_0, \zeta) + \sqrt{n} h' \mathbb{P}_n \tilde{l}_\beta(\theta_0(\zeta)) - \tfrac{1}{2} h' \tilde{I}_0 h + o_{P_0}^{\overline{\Xi}}(1) \\
&= \tfrac{1}{2} \overline{\beta}_n' \tilde{I}_0 \overline{\beta}_n - \tfrac{1}{2}(\overline{\beta}_n - h)' \tilde{I}_0(\overline{\beta}_n - h) + o_{P_0}^{\overline{\Xi}}(1),
\end{aligned}
$$

therefore,

$$
\begin{aligned}
PLR_n(M) &= \iint_{\|h\| \leq M} (pl_n(\beta_0 + h/\sqrt{n}, \zeta) - pl_n(\beta_0, \zeta)) \, dQ_\zeta(h) \, dJ(\zeta) \\
&= \iint_{\|h\| \leq M} \exp(\tfrac{1}{2} \overline{\beta}_n' \tilde{I}_0 \overline{\beta}_n \\
&\qquad\qquad - \tfrac{1}{2}(\overline{\beta}_n - h)' \tilde{I}_0(\overline{\beta}_n - h) + o_p^{\overline{\Xi}}(1)) \, dQ_\zeta(h) \, dJ(\zeta) \\
&= \overline{PLR}_n(M) + o_p(1)
\end{aligned}
$$

where the last equality follows from $\overline{PLR}_n(M) = O_p(1)$, by arguments analogous to those used in (16) above. The proof for Part (1) is now completed.

For Part (2), since $h \sim Q_\zeta = N(0, c\tilde{I}_0^{-1})$,

$$\overline{PLR}_n = \int_{\zeta \in \Xi} \xi_n(\zeta) \, dJ(\zeta),$$

with

$$\xi_n(\zeta) = \int \exp\left(\frac{1}{2}\overline{\beta}_n' \tilde{I}_0 \overline{\beta}_n - \frac{1}{2}(\overline{\beta}_n - h)' \tilde{I}_0 (\overline{\beta}_n - h)\right) dQ_\zeta(h)$$

$$= (2\pi)^{-p/2} \det^{1/2}(\tilde{I}_0/c)$$

$$\times \int \exp\left[\frac{1}{2}\left\{\overline{\beta}_n' \tilde{I}_0 \overline{\beta}_n - (h - \overline{\beta}_n)' \tilde{I}_0 (h - \overline{\beta}_n) - \frac{h' \tilde{I}_0 h}{c}\right\}\right] dh$$

$$= (1+c)^{-p/2} \exp\left(\frac{1}{2}\frac{c}{1+c}\overline{\beta}_n' \tilde{I}_0 \overline{\beta}_n\right),$$

where the last equality holds by integrating out a normal density.

For Part (3), it follows from Lemma 2 and assumption B3 that $\sqrt{n}\|\hat{\beta}_n(\zeta) - \beta_0\| = O_{\overline{P}_0}^{\overline{\Xi}}(1)$, and reapplication of Lemma 2 and the arg max theorem yields $\sqrt{n}(\hat{\beta}_n(\zeta) - \beta_0) = \tilde{I}_\beta(\theta_0(\zeta))^{-1}\sqrt{n}\mathbb{P}_n\tilde{I}_\beta(\theta_0(\zeta)) + o_{\overline{P}_0}^{\overline{\Xi}}(1)$. Part (3) now follows.

For the proof of Part (4) and Part (5), it suffices to show that $W_n(\zeta) - R_n(\zeta) = o_{\overline{P}_0}^{\overline{\Xi}}(1)$ and $R_n(\zeta) - LR_n(\zeta) = o_{\overline{P}_0}^{\overline{\Xi}}(1)$. These results follow from Donsker properties and standard arguments. We omit the details. The proof of Lemma 5 is thus complete. $\square$

PROOF OF COROLLARY 1.    The proof is similar to the proof of Theorem 1. We omit the details. $\square$

PROOF OF COROLLARY 2.    The proof follows the same lines as the proof of Part (2)(iii) of Lemma 5, with

$$W(q_\zeta, \zeta) = (2\pi)^{-p/2} \det^{1/2}\left(\frac{1+c}{c}\tilde{I}_0\right)$$

$$(17) \qquad \times \int \exp\left[-\frac{1+c}{2c}\left(\lambda - \frac{c}{1+c}\overline{\beta}_n\right)' \tilde{I}_0 \left(\lambda - \frac{c}{1+c}\overline{\beta}_n\right)\right.$$

$$\left. - \lambda'\langle q_\zeta - \tilde{q}_\zeta, P\tilde{i}_\eta^\star \tilde{i}_\eta (q_\zeta - \tilde{q}_\zeta)'\rangle_\eta \lambda\right] d\lambda,$$

where det is the determinant of a matrix, $\langle \cdot, \cdot \rangle_\eta$ is the inner product defined on $\overline{\mathcal{H}}_\eta$, and $W(q_\zeta, \zeta) \leq 1$ since $\langle q_\zeta - \tilde{q}_\zeta, P\tilde{i}_\eta^\star \tilde{i}_\eta (q_\zeta - \tilde{q}_\zeta)'\rangle_\eta$ is nonnegative definite. $\square$

LEMMA 6.   *Under assumptions* A–C, *the densities* $\ell_n(\overline{\psi}_0 + h/\sqrt{n}, \zeta)$ *and* $\int \ell_n(\overline{\psi}_0 + h/\sqrt{n}, \zeta) \, dQ_\zeta(h) \, dJ(\zeta)$ *are contiguous to the densities* $l_n^0$. *As a consequence, the results of Lemma* 5 *still hold under local alternatives* $\{P_{\overline{\psi}_0 + h/\sqrt{n}, \zeta}\}$ *and* $\{\int P_{\overline{\psi}_0 + h/\sqrt{n}, \zeta} \, dQ_\zeta(h) \, dJ(\zeta)\}$.

PROOF.   Assumption C implies that a LAN (local asymptotic normal) expansion for the log-likelihood ratio holds immediately by Lemma 3.10.11 of van der Vaart and Wellner (1996):

$$\Lambda_{n\zeta} \equiv \log\left(\frac{dP_{\psi_0 + h/\sqrt{n}, \zeta}^n}{dP_{\psi_0, \zeta}^n}\right) = \frac{1}{\sqrt{n}} \sum_{i=1}^n A_\zeta h(X_i) - \frac{1}{2}\|A_\zeta h\|^2 + o_{P_0}(1).$$

It follows from LAN that $\Lambda_{n\zeta} \to_d W_\zeta$, where $W_\zeta \sim N(-1/2\|A_\zeta h\|^2, \|A_\zeta h\|^2)$, under $P_0$. Therefore, under $P_0$,

$$\exp(\Lambda_{n\zeta}) \equiv \frac{dP_{\overline{\psi}_0 + h/\sqrt{n}, \zeta}^n}{dP_0^n} \xrightarrow{d} \exp W_\zeta.$$

$P_0(\exp(W_\zeta)) = 1$, using the formula for the moment generating function of the normal distribution. By Le Cam's first lemma [van der Vaart (1996), page 88], we conclude that the sequences of probability measures $\{P_{\overline{\psi}_0 + h/\sqrt{n}, \zeta}\}$ and $\{P_0\}$ are contiguous, for every $\zeta \in \Xi$. Consequently the convergence in probability that holds under $P_0$ also holds under $\{P_{\overline{\psi}_0 + h/\sqrt{n}, \zeta}\}$ and vice versa. Similarly, since $P(e\chi) = 1$ using the formula for the moment generating function of the $\chi^2$ distribution, we conclude that the sequences $\{\int P_{\overline{\psi}_0 + h/\sqrt{n}, \zeta}^n \, dQ_\zeta(h) \, dJ(\zeta)\}$ and $P_0^n$ are contiguous.   □

PROOF OF THEOREM 2.   We define a $\sqrt{n}$-neighborhood of $\beta_0$ as a collection of sequences $\beta_n(h_\beta) = \beta_0 + h_\beta/\sqrt{n} + o(n^{-1/2})$, for $h_\beta \in \mathbb{R}^p$. A $\sqrt{n}$ neighborhood of $\eta$ is similarly defined as $\eta_n(h_\eta) = \eta + h_\eta/\sqrt{n} + o(n^{-1/2})$, for $h_\eta \in \mathcal{H}_\eta$. With a minor abuse of notation, a local form of the hypotheses can be written as:

$$(18) \qquad H_0 : \overline{\psi} = \overline{\psi}_0 \quad \text{vs.} \quad H_1 : \overline{\psi} = \overline{\psi}_0 + h_1/\sqrt{n},$$

where $h_1 \in \mathbb{R}^p \times \mathcal{H}_\eta$ takes the value $(h_{\beta 1}, h_{\eta 1})$, with $h_{\eta 1} = \tilde{q}_\zeta' h_{\beta 1}$. We note that the least favorable direction $\tilde{q}_\zeta$ is invariant under the choice of $\phi_\zeta$, and, as a consequence, the contiguous alternative $H_1$ is also invariant under the choice of $\phi_\zeta$.

Define

$$(19) \qquad LR_n \equiv \frac{\int \ell_n(\overline{\psi}_0 + h_1/\sqrt{n}, \zeta) \, dQ_\zeta(h_{\beta 1}) \, dJ(\zeta)}{\ell_n^0}.$$

A test defined by $LR_n$ is

$$\tilde{\gamma}_n = \begin{cases} 1, & \text{if } LR_n > \tilde{k}_{\alpha n}, \\ \tilde{\lambda}_n, & \text{if } LR_n = \tilde{k}_{\alpha n}, \\ 0, & \text{if } LR_n < \tilde{k}_{\alpha n}, \end{cases}$$

where $\tilde{k}_{\alpha n} > 0$, $\tilde{\lambda}_n \in [0,1]$ are constants such that the rejection probability is $\alpha$ under the null. For notational simplicity, let $P_1^n = \int P_{\overline{\psi}_0 + h_1/\sqrt{n}, \zeta}^n \, dQ_\zeta(h_{\beta 1}) \, dJ(\zeta)$. By the Neyman–Pearson lemma, for all $n \geq 1$ and any test $\phi_n$ with level $\alpha$, with a minor abuse of notation,

$$\lim_{n \to \infty} \int \phi_n \left\{ \int \ell_n(\overline{\psi}_0 + h_1/\sqrt{n}, \zeta) \, dQ_\zeta(h_{\beta 1}) dJ(\zeta) \right\} dP_1^n \tag{20}$$

$$\leq \lim_{n \to \infty} \int \tilde{\gamma}_n \left\{ \int \ell_n(\overline{\psi}_0 + h_1/\sqrt{n}, \zeta) \, dQ_\zeta(h_{\beta 1}) \, dJ(\zeta) \right\} dP_1^n$$

$$= \lim_{n \to \infty} \int I(LR_n > \tilde{k}_{\alpha n}) \tag{21}$$

$$\times \left\{ \int \ell_n(\overline{\psi}_0 + h_1/\sqrt{n}, \zeta) \, dQ_\zeta(h_{\beta 1}) \, dJ(\zeta) \right\} dP_1^n$$

$$= \lim_{n \to \infty} \int \left\{ \int I(PLR_n > \tilde{k}_{\alpha n}) \, dP_{\overline{\psi}_0 + h_1/\sqrt{n}, \zeta}^n \right\} dQ_\zeta(h_{\beta 1}) \, dJ(\zeta) \tag{22}$$

$$= \lim_{n \to \infty} \int \left\{ \int I(EW_n > \tilde{k}_{\alpha n}) \, dP_{\overline{\psi}_0 + h_1/\sqrt{n}, \zeta}^n \right\} dQ_\zeta(h_{\beta 1}) \, dJ(\zeta), \tag{23}$$

where (21) follows since $LR_n$ has an absolutely continuous asymptotic distribution under the contiguous alternative $H_1$ and by Fubini's theorem. (22) follows since $PLR_n - LR_n = o_P(1)$ under $H_1$, which will be established at the end of the proof. (23) follows from Lemma 6. The results for $ER_n$ and $ELR_n$ also follow from Lemma 6. By Fubini's theorem, we obtain $\limsup_{n \to \infty} \int \{ \phi_n (P_{\overline{\psi}_0 + h_1/\sqrt{n}, \zeta}^n) \} dQ_\zeta(h_{\beta 1}) \, dJ(\zeta) \leq \lim_{n \to \infty} \int \{ \int I(EW_n > \tilde{k}_{\alpha n}) \, dP_{\overline{\psi}_0 + h_1/\sqrt{n}, \zeta}^n \} dQ_\zeta(h_{\beta 1}) \, dJ(\zeta)$, which implies that the proposed tests have the greatest weighted average power asymptotically in the class of all tests of asymptotic significance level $\alpha$, against the alternative $P_{\overline{\psi}_0 + h/\sqrt{n}, \zeta}^n$.

To show $PLR_n - LR_n = o_P(1)$ under $H_1$, it suffices to show $PLR_n - LR_n = o_P(1)$ under the null by Lemma 6. Define $LR_n(M) \equiv \int_{\zeta \in \Xi} \int_{\|h\| \leq M} \ell_n(\overline{\psi}_0 + h_1/\sqrt{n}, \zeta) \, dQ_\zeta(h) \, dJ(\zeta) / \ell_n(\overline{\psi}_0, \zeta)$, and note that $\forall M : 0 < M < \infty$, $|PLR_n - LR_n| \leq |PLR_n - PLR_n(M)| + |PLR_n(M) - LR_n(M)| + |LR_n - LR_n(M)|$. Hence it suffices to show that: (i) $|PLR_n - PLR_n(M)| \to_{P_0} 0$, (ii) $|\overline{LR}_n - \overline{LR}_n(M)| \to_{P_0} 0$ and (iii) $|PLR_n(M) - LR_n(M)| \to_{P_0} 0$, as $n \to \infty$. Part (i) was shown in Lemma 5. Part (ii) can be similarly established by taking $M$ large enough and using assumption A.

To show Part (iii), we take Taylor expansion of $\log \ell_n(\overline{\psi}_0 + h_1/\sqrt{n}, \zeta)$ at $(\overline{\psi}_0, \zeta)$ with respect to $h_\beta$ along the direction $\tilde{q}_\zeta$, which leads to the following expansion in the least favorable submodel:

$$\log \ell_n\left(\overline{\psi}_0 + \frac{h_1}{\sqrt{n}}, \zeta\right) = \log \ell_n(\overline{\psi}_0, \zeta) + \sqrt{n} h'_{\beta 1} \mathbb{P}_n \dot{\ell}(\beta_0, \overline{\psi}_0, \zeta)$$

$$+ \tfrac{1}{2} h'_{\beta 1} \mathbb{P}_n \ddot{\ell}(\tilde{\beta}, \tilde{\psi}, \zeta) h_{\beta 1}.$$

On the right-hand side, we can replace $\mathbb{P}_n \dot{\ell}(\beta_0, \overline{\psi}_0, \zeta)$ by $\mathbb{P}_n \dot{\ell}_\beta(\overline{\psi}_0, \zeta) + o_{P_0}^{\Xi}(1)$, and $\mathbb{P}_n \ddot{\ell}(\tilde{\beta}, \tilde{\psi}, \zeta)$ by $-\tilde{I}_\beta(\overline{\psi}_0, \zeta) + o_{P_0}^{\Xi}(1)$, according to assumption B2. Comparing the above display and Lemma 2 with $\tilde{\beta}_n \equiv h_{\beta 1}/\sqrt{n}$, we obtain Part (iii). $\square$

PROOF OF THEOREM 3. The equivalence of the three tests under local alternatives is shown in Lemma 6. To show their asymptotic distribution, a key step is to establish that $\overline{\beta}_n$ converges under $P_{\overline{\psi}_0 + h/\sqrt{n}, \zeta_1}^n$ in distribution to the process $\zeta \mapsto \mathbb{G}(\theta_0(\zeta)) + \nu_\star(h_\beta, \zeta, \zeta_1)$, where $\nu_\star(h_\beta, \zeta, \zeta_1) \equiv P_0 \tilde{l}_\beta(\theta_0(\zeta)) \tilde{l}_\beta(\theta_0(\zeta_1))' h_\beta$, by Theorem 3.10.12 in van der Vaart and Wellner (1996). The result follows by Lemma 6 and the continuous mapping theorem. $\square$

PROOF OF THEOREM 4. The equivalence of the three tests under local alternatives is shown in Lemma 6. Since the sequences of densities $\int \ell_n(\overline{\psi}_0 + h/\sqrt{n}, \zeta) \, dQ_\zeta(h) \, dJ(\zeta)$ are contiguous to the density $l_n^0$, we have

$$\left( ELR_n, \frac{\int dP_{\overline{\psi}_0 + h/\sqrt{n}, \zeta}^n \, dQ_\zeta(h) \, dJ(\zeta)}{dP_0^n} \right) \rightsquigarrow_d (e\chi(c), e\chi(c)),$$

under $P_0$. Then $ELR_n \to_d r\chi(c)$ under $\int dP_{\overline{\psi}_0 + h/\sqrt{n}, \zeta}^n \, dQ_\zeta(h) \, dJ(\zeta)$, by Le Cam's third lemma. $\square$

PROOF OF LEMMA 3. The proof mainly involves an argument that for an arbitrary, possibly random sequence $\{\zeta_n\}$, the distance between the minimizer of the Kullback–Leibler information and $\hat{\theta}_n(\zeta_n)$ goes to zero. Consequently, the assertion of Lemma 3 follows from the arbitrariness of the sequence $\zeta_n$ and Slutsky's theorem. We omit the details. $\square$

PROOF OF LEMMA 4. The proof mainly involves a uniform "peeling device" with an adaptation of the proof of Theorem 3.2 given in Murphy and van der Vaart (1999), which details we omit. $\square$

LEMMA 7. (1) *Assume* $\phi_\zeta : \mathbb{D}_\phi \subset \mathbb{D} \mapsto \mathbb{E}_\psi \subset \mathbb{E}$ *is one-to-one, continuously invertible and onto and* $\psi_\zeta : \mathbb{E}_\psi \subset \mathbb{E} \mapsto \mathbb{F}$ *is one-to-one, continuously invertible and onto, then* $\psi_\zeta \circ \phi_\zeta : \mathbb{D}_{\phi_\zeta} \mapsto \mathbb{F}$ *is one-to-one, continuously invertible and onto.* (2) *Assume* $\phi_\zeta : \mathbb{D}_\phi \subset \mathbb{D} \mapsto \mathbb{E}_\psi \subset \mathbb{E}$ *is uniformly Fréchet differentiable at* $\theta \in \mathbb{D}_\psi$ *and* $\psi_\zeta : \mathbb{E}_\psi \subset \mathbb{E} \mapsto \mathbb{F}$ *is uniformly Fréchet differentiable at* $\phi_\zeta(\theta)$ *over* $\zeta \in \Xi$. *Then* $\psi_\zeta \circ \phi_\zeta : \mathbb{D}_\phi \mapsto \mathbb{F}$ *is uniformly Fréchet differentiable at* $\theta$ *with derivative* $\psi'_\zeta(\phi_\zeta(\theta)) \circ \phi'_\zeta(\theta)$.

PROOF. For Part (1), it suffices to note that $\|\dot\psi_\zeta(\phi_\zeta(\theta))(\dot\phi_\zeta(\theta)(h))\| \geq c_1 \|\dot\phi_\zeta(\theta)(h)\| \geq c_1 c_2 \|h\|$. For Part (2), we note that, $\psi_\zeta \circ \phi_\zeta(\theta + th) - \psi_\zeta \circ \phi_\zeta(\theta) = \psi_\zeta(\phi_\zeta(\theta) + tk_t) - \psi_\zeta(\phi_\zeta(\theta))$, where $k_t = \{\phi_\zeta(\theta + th) - \phi_\zeta(\theta)\}/t$. So we rewrite the uniform Fréchet difference as $\psi_\zeta(\phi_\zeta(\theta + h))(\cdot) - \psi_\zeta(\phi_\zeta(\theta))(\cdot) = \dot\psi_\zeta(\phi_\zeta(\theta))(\phi_\zeta(\theta + h) - \phi_\zeta(\theta)) + o^\Xi(\|\phi_\zeta(\theta + h) - \phi_\zeta(\theta)\|) = \dot\psi_\zeta(\phi_\zeta(\theta)) \times \dot\phi_\zeta(\theta)(h) + o^\Xi(\|h\|)$. □

LEMMA 8. *Let* $A_\zeta = T_\zeta + K_\zeta : \mathbb{D} \mapsto \mathbb{E}$ *be a linear operator between Banach spaces, where* $T_\zeta$ *is onto and there exists* $c_1 > 0$, *such that* $\|T_\zeta h\| \geq c_1 \|h\|$ *for all* $h \in \mathbb{D}$ *and* $\zeta \in \Xi$, *and* $K_\zeta$ *is uniformly compact, that is,* $\bigcup_{\zeta \in \Xi} \bigcup_{\|h\| \leq 1} K_\zeta h$ *is compact. Then if* $N(A_\zeta) = \{0\}$ *for all* $\zeta \in \Xi$, *then* $A_\zeta$ *is onto and there exists* $c_2 > 0$ *such that* $\|A_\zeta h\| \geq c_2 \|h\|$, $\forall \zeta \in \Xi$ *and all* $h \in \mathbb{D}$.

PROOF. Since, for an arbitrary random sequence $\zeta_n$, $T_{\zeta_n}^{-1}$ is continuous, the operator $T_{\zeta_n}^{-1} K : \mathbb{E} \mapsto \mathbb{D}$ is compact. Hence $I + T_{\zeta_n}^{-1} K_{\zeta_n}$ is one-to-one and therefore also onto be a result of Riesz for compact operators. Thus $T_{\zeta_n} + K_{\zeta_n}$ is also onto. We will be done if we can show that $I + T_{\zeta_n}^{-1} K_{\zeta_n}$ is continuously invertible, since that would imply that $(T_{\zeta_n} + K_{\zeta_n})^{-1} = (I + T_{\zeta_n}^{-1} K_{\zeta_n})^{-1} T_{\zeta_n}^{-1}$ is bounded. The remainder of the proof follows the proof of Lemma 6.17 in [Kosorok (2008)]. □

LEMMA 9. *Suppose that* $U_n(\psi, \zeta)(h) = \mathbb{P}_n \nu(\psi, \zeta)(h)$ *and* $U(\psi, \zeta)(h) = P\nu(\psi, \zeta)(h)$ *for given* $P$-*measurable functions* $\nu(\psi, \zeta)(h)$ *indexed by* $\Psi \times \Xi$ *and an arbitrary index set* $\mathcal{H}_\eta$. *Assume* $\psi = \psi_0$, $\nu(\psi_0, \zeta)(h) = \nu(\psi_0)(h)$. *If* $\hat\psi_n(\zeta) = \psi_0 + o_P^\Xi(1)$, *the class of functions* $\{\nu(\psi, \zeta)(h) - \nu(\psi_0)(h) : \|\psi - \psi_0\| < \delta, h \in \mathcal{H}_\eta, \zeta \in \Xi\}$ *is* $P$-*Donsker for some* $\delta > 0$ *and* $\sup_{\zeta \in \Xi, h \in \mathcal{H}_\eta} P_0(\nu(\psi, \zeta)(h) - \nu(\psi_0)(h))^2 \to 0$, *as* $\psi \to \psi_0$, *then* $\sup_{\zeta \in \Xi} \|\sqrt{n}(U_n - U)(\hat\psi_n(\zeta), \zeta) - \sqrt{n}(U_n - U)(\psi_0, \zeta)\| = o_P^\Xi(1 + \sqrt{n}\|\hat\psi_n(\zeta) - \psi_0\|)$.

PROOF. This is a "uniform" version of Lemma 3.3.5 in van der Vaart and Wellner [(1996)]. Let $\Psi_\delta \equiv \{\psi : \|\psi - \psi_0\| < \delta\}$ and define an extraction function $f : \ell^\infty(\Psi_\delta \times \Xi \times \mathcal{H}_\eta) \times \Psi_\delta \mapsto \ell^\infty(\mathcal{H}_\eta \times \Xi)$ as $f(z, \psi, \zeta)(h) \equiv z(\psi, \zeta, h)$, where $z \in \ell^\infty(\Psi_\delta \times \mathcal{H}_\eta \times \Xi)$. Since $f$ is continuous at every point $(z, \psi_1, \zeta)$, we have $\sup_{h \in \mathcal{H}_n, \zeta \in \Xi} |z(\psi, \zeta, h) - z(\psi_1, \zeta, h)| \to 0$ as $\psi \to \psi_1$. Define the stochastic process $Z_n(\psi, \zeta, h) \equiv \mathbb{G}_n(\nu(\psi, \zeta)(h) - \nu(\psi_0, \zeta)(h))$ indexed by $\Psi_\delta \times \Xi \times \mathcal{H}_\eta$. By assumptions, $Z_n$ converges weakly in $\ell^\infty(\Psi_\delta \times \Xi \times \mathcal{H}_\eta)$ to a tight Gaussian process $Z_0$ with continuous sample paths with respect to the metric $\rho_\zeta$ defined by $\rho_\zeta^2((\psi_1, \zeta, h_1), (\psi_2, \zeta, h_2)) = P(\nu(\psi_1, \zeta)(h_1) - \nu(\psi_0, \zeta)(h_1) - \nu(\psi_2, \zeta)(h_2) + \nu(\theta_0, \zeta)(h_2))^2$, at fixed $\zeta$. Since as assumed, $\sup_{h \in \mathcal{H}_n, \zeta \in \Xi} \rho_\zeta((\psi, h), (\psi_0, h)) \to 0$, we have that $f$ is continuous at almost all sample paths of $Z_0$ uniformly over $\zeta \in \Xi$. The result now follows by Slutsky's theorem and the continuous mapping theorem. $\square$

**Acknowledgments.** We thank the referees and editors for several helpful comments that led to an improved paper.

R. SONG
J. P. FINE
DEPARTMENT OF BIOSTATISTICS
UNIVERSITY OF NORTH CAROLINA
    AT CHAPEL HILL
CHAPEL HILL, NORTH CAROLINA 27599-7420
USA
E-MAIL: rsong@bios.unc.edu
          jfine@bios.unc.edu

M. R. KOSOROK
DEPARTMENTS OF BIOSTATISTICS
    AND STATISTICS & OPERATIONS RESEARCH
UNIVERSITY OF NORTH CAROLINA
    AT CHAPEL HILL
CHAPEL HILL, NORTH CAROLINA 27599-7420
USA
E-MAIL: kosorok@unc.edu